\documentclass[10pt,twocolumn,twoside]{IEEEtran} 
\usepackage[utf8]{inputenc}

\IEEEoverridecommandlockouts

\usepackage{graphicx}
\usepackage{amsmath}
\usepackage{amssymb}
\usepackage{amsfonts}
\usepackage{mathrsfs}
\usepackage{dsfont}
\makeatletter
\def\amsbb{\use@mathgroup \M@U \symAMSb}
\makeatother
\usepackage{color}

\usepackage[usenames,dvipsnames]{xcolor}
\usepackage{epstopdf}

\newtheorem{theorem}{Theorem}
\newtheorem{remark}{Remark}
\newtheorem{definition}{Definition}
\newtheorem{lemma}{Lemma}
\newtheorem{proposition}{Proposition}
\newtheorem{corollary}{Corollary}
\newtheorem{assumption}{Assumption}
\newtheorem{design condition}{Design condition}
\newtheorem{problem}{Problem}
 
\newcommand{\norm}[1]{\left\lVert#1\right\rVert}
\usepackage{enumerate}
\usepackage{balance}

\DeclareMathOperator*{\argmin}{arg\,min}

\allowdisplaybreaks

\usepackage{mathtools}
\DeclarePairedDelimiter\floor{\lfloor}{\rfloor}

\newcommand{\Mod}[1]{\mathrm{mod} #1}

\newcommand{\ak}[1]{{\color{black} #1}}
\newcommand{\an}[1]{{\color{black} #1}}
\newcommand{\ank}[1]{{\color{black} #1}}
\newcommand{\ka}[1]{{\color{black} #1}}

\title{On the stability properties of power networks with time-varying inertia}

\author{Andreas Kasis,\thanks{This work was funded by the European Union’s Horizon 2020 research and innovation program under grant agreements No. 891101 (SmarTher Grid) and No. 739551 (KIOS CoE), and from the Republic of Cyprus through the Directorate General for European Programs, Coordination, and Development.}
 Stelios Timotheou and  Marios Polycarpou
 \thanks{Andreas Kasis, Stelios Timotheou and 
Marios Polycarpou are with the KIOS Research and Innovation Center of Excellence and the Department of Electrical and Computer Engineering, University of Cyprus, Cyprus; e-mails: \{kasis.andreas, timotheou.stelios, mpolycar\}@ucy.ac.cy.}
\thanks{A preliminary version of this work has appeared in \cite{CDC_2021}. This manuscript provides additional results and the analytic proofs of the main results. Moreover, it includes additional discussion and simulations that demonstrate
the impact of the proposed analysis.}
}

\begin{document}

\maketitle

%
%
%
%
%
%
%
%
%
%
%

%

\begin{abstract}

A major transition in modern power systems is the replacement of conventional generation units with renewable sources of energy. 
The latter results in lower rotational  inertia which compromises the stability of the power system, as testified by the growing number of frequency incidents.
To resolve this problem, numerous studies have proposed the use of virtual inertia to improve the stability properties of the power grid.
In this study, we consider how inertia variations, resulting from the application of control action associated with virtual inertia, may affect the stability properties of the power network within the primary frequency control timeframe.
We consider the interaction between the frequency dynamics and a broad class of power supply dynamics in the presence of time-varying inertia  and provide locally verifiable conditions, that enable scalable designs, such that stability is guaranteed.
To complement the presented stability analysis and highlight the dangers arising from varying inertia,  we provide analytic conditions that enable to deduce instability from single-bus inertia fluctuations. 
Our analytical results are validated with simulations on the Northeast Power Coordinating Council (NPCC) 140-bus and the IEEE New York/New England 68-bus systems, where  we demonstrate  how  inertia variations may induce large frequency oscillations and show that the application of the proposed conditions yields a stable response.
\end{abstract}


\section{Introduction}

\textbf{Motivation and literature review:} 
The electric power grid is currently undergoing through a major transformation due to the growing penetration of renewable sources of energy  \cite{lund2006large}, \cite{ipakchi2009grid}.
As a result, conventional bulk generation is expected to be slowly replaced by renewable generation.
However, retiring synchronous generation  lowers the rotational inertia of the power system, which has been a key reason for the power grid's stability over the years \cite{tielens2016relevance}.
In addition, renewable generation is intermittent, causing more  frequent  generation-demand imbalances that may harm the power quality and even cause blackouts \cite{lalor2004dynamic}. 
Hence, novel challenges are introduced towards enabling the stable and robust operation of the power grid.




The inertia of the power system represents its capability to store and inject kinetic energy, serving as an energy buffer that can slow down the frequency dynamics.
The latter aids to  effectively avoid undesirable incidents, such as excessive load-shedding or large-scale blackouts.
A low power system inertia is associated with larger frequency deviations following a disturbance event, such as loss of generation or tie line faults \cite{ulbig2014impact}.
The degrading effects of low inertia levels on the power system  stability have already been reported by system operators \cite{entso_report}.
To mitigate these effects, several studies proposed the introduction of  virtual inertia   in the power grid, i.e. schemes that aim to resemble the inertial response of machines by injecting power in proportion to the rate of change of frequency.
In particular, \cite{fang2017distributed} proposed the use of the energy stored in the DC-link capacitors of grid-connected power converters to emulate the inertia of the power system.
Control schemes of varying complexity that aim to make converters behave in a similar way as synchronous
machines  and provide virtual inertia to the power grid have been investigated in \cite{zhong2010synchronverters}, \cite{arghir2018grid}, \cite{van2010virtual}, \cite{torres2009virtual}, \cite{alipoor2014power}.
These schemes require some form of energy storage, such as batteries, to provide the required  energy. 
In addition, \cite{chen2017integration} proposed the integration of DC microgrids  into the power network as virtual synchronous machines (VSMs).
Moreover, \cite{lopes2014self} designed a self tuning mechanism to optimize the operation of VSMs. 
The scheme proposed in  \cite{morren2006wind}  uses the kinetic energy stored in the rotating mass of wind turbine blades to emulate inertia and provide frequency support.
Furthermore, \cite{hughes2005control}, \cite{ekanayake2004comparison} demonstrated how power converters can mimic the inertial response of a synchronous machine for  wind turbines interconnected to the grid via doubly fed induction generators.
The optimal placement of virtual inertia is  considered in \cite{poolla2017optimal}.
 Comprehensive reviews on virtual inertia and virtual synchronous machine schemes are available in \cite{tamrakar2017virtual} and \cite{d2013virtual} respectively.


An open research question  that requires further attention concerns the effect of varying inertia on the stability properties of the power network. 
Variations in  inertia may arise due to 
control action on virtual inertia.
Controllable virtual inertia,  possibly coupled with the  power grid dynamics,  may induce additional challenges for the stability of the power grid.
In addition, it may pose security threats
if inertia
is maliciously controlled
   to destabilize the power grid.
The time-varying nature of inertia has been pointed out and studied in \cite{ulbig2014impact}.
Moreover, \cite{misyris2018robust} considered the robustness properties of power networks with time-varying inertia and frequency damping, while  \cite{hidalgo2018frequency} considered a hybrid model of the power network with time-varying inertia and applied model predictive control approaches to optimize the power inputs. 
In \cite{borsche2015effects}, the authors considered the effect of inertia variations in the frequency response.
Stabilizing controllers under various inertia modes are designed in \cite{srivastava2022learning} using a linear quadratic regulator approach.
Furthermore, \cite{wu2018frequency} and \cite{alipoor2013distributed} proposed time-varying control gains associated with the inertia and damping of wind turbines and the inertia of virtual synchronous generators respectively.
However, a systematic approach on how time-varying inertia affects the stability properties of power networks is currently missing from the literature.
  This may lead to suitable guidelines for the design of  virtual inertia schemes that  offer improved stability and security properties.

\textbf{Contribution:} 
This study investigates the impact of time-varying inertia on the behaviour and stability properties of the power network within the primary frequency control timeframe. 
In particular, we consider the interaction between the frequency dynamics and nonlinear power supply dynamics in the presence of time-varying inertia.
We  study the solutions of this system and analytically deduce local asymptotic stability  under two proposed conditions.
The first condition, inspired from \cite{kasis2016primary}, requires that the aggregate power supply dynamics at each bus satisfy a passivity related property.
The second condition sets a bound on the maximum rate of growth of inertia that depends on the local power supply dynamics.
Both conditions are locally verifiable and applicable to general network topologies and 
 enable practical   guidelines for the design of virtual inertia control schemes that enhance the reliability and security of the power network.

To demonstrate the applicability and improve the intuition on the proposed analysis, we provide examples of power supply dynamics that fit within the presented framework and explain how the maximum allowable rate of inertia variation may be deduced for these cases.
In addition, when linear power supply dynamics are considered, we show how the conditions can be efficiently verified by solving a suitable linear matrix inequality optimization problem.

Our stability results are complemented by analytical results that demonstrate how single-bus inertia variations may induce unstable behaviour. 
In particular, we provide conditions such  that for any, arbitrarily small, deviation from the equilibrium frequency at a single bus, there exist local inertia trajectories that result in substantial deviations in the power network frequency. 
The latter coincides with the definition of instability (e.g. \cite[Dfn. 4.1]{haddad2011nonlinear}).
\an{These results confirm that inertia-induced instability is a fundamental, network-independent phenomenon rooted in frequency dynamics, offering insight into its underlying mechanisms.}

%
%



Numerical simulations on the NPCC 140-bus and the IEEE New York/New England 68-bus systems demonstrate the potentially destabilizing effects of varying inertia and validate our analytical results in a realistic setting.
 In particular,  we demonstrate how varying virtual inertia, at a single or multiple buses, may induce large frequency oscillations and compromise the stability of the power network.
The latter provides further motivation for the study and regulation of virtual inertia schemes.
In addition, we present how the application of the proposed conditions yields a stable response.

To the authors best knowledge, this is the first work that:
\begin{enumerate}[(i)]
\item Analytically studies the behaviour of power networks under varying inertia and proposes decentralized conditions on the local power supply dynamics and inertia trajectories such that stability is guaranteed.
\item Demonstrates, under broadly applicable conditions, that single-bus inertia variations may induce large frequency fluctuations from any, arbitrarily small, frequency deviation from the equilibrium. 
Such behaviour is characterized as unstable.
\end{enumerate}

\textbf{Paper structure:} 
In Section \ref{sec:Problem_formulation} we present the power network model, a general description for the power supply dynamics and a statement of the problem that we consider.
In Section \ref{sec:conditions} we present our proposed conditions on the power supply dynamics.
Section \ref{Sec: Convergence} contains the conditions on virtual inertia trajectories and the main result associated with the stability of the power network.
In addition, in Section \ref{Sec: Discussion} we provide additional intuition on the stability result and discuss various application examples.
Moreover, in Section \ref{Sec_Instability} we present our inertia induced instability analysis.
Our analytical results are verified with numerical simulations in Section \ref{Simulation_NPCC} and conclusions are drawn in Section \ref{Sec: Conclusions}.
Proofs of the main results are provided in the Appendix.

\textbf{Notation:}
Real, positive real, integer and positive natural numbers are denoted by $\mathbb{R}, \mathbb{R}_+, \mathbb{Z}$ and $\mathbb{N}_+$ respectively.
The set of n-dimensional vectors with real entries is denoted by $\mathbb{R}^n$.
For $a \in \mathbb{R}, b \in \mathbb{R} \setminus \{0\}$, $a$ modulo $b$ is denoted by~$\Mod(a, b)$ and defined as
$\Mod(a, b) = a - b \floor{\frac{a}{b}}
$, where for $x \in \mathbb{R}$, $\floor{x} = \sup \{m \in \mathbb{Z}: m \leq x\}$.
The $p$-norm of a vector $x \in \mathbb{R}^n$ is given by $\norm{x}_p = (|x_1|^p + \dots + |x_n|^p)^{1/p}, 1 \leq p < \infty$.
A function $f: \mathbb{R}^n \rightarrow \mathbb{R}^m$ is said to be locally Lipschitz continuous at $x$ if there exists some neighbourhood $X$ of $x$ and some constant $L$ such that $\norm{f(x) - f(y)} \leq L \norm{x - y}$ for all $y \in X$, where $\norm{.}$ denotes any $p$-norm, and globally Lipschitz if the condition holds on the whole domain of $y$.
The Laplace transformation of a signal $x(t), x : \mathbb{R} \rightarrow \mathbb{R}$, is denoted
by $\hat{x}(s) = \int_0^{\infty} x(t) e^{-st} dt$.
A function $f: \mathbb{R}^n \rightarrow \mathbb{R}$ is called positive semidefinite if $f(x) \geq 0$ for all $x \in \mathbb{R}^n$.
The function \textbf{sin}$(x)$ gives the sinusoid of each element in $x \in \mathbb{R}^n$, i.e. for $x = [x_1 \; x_2 \dots x_n]^T$, \textbf{sin}$(x) = [\sin (x_1) \; \sin (x_2) \dots \sin (x_n)]^T$.
For a matrix $A \in \mathbb{R}^{n \times p}$, $A_{kl}$ corresponds to the element in the $k$th row and $l$th column of $A$.
A matrix $A \in \mathbb{R}^{n \times n}$ is called diagonal if $A_{ij} = 0$ for all $i \neq j$ and positive (negative) semi-definite, symbolized with $A \succeq 0$ (respectively $A \preceq 0)$, if $x^T A x \geq  0$ (respectively $x^T A x \leq  0$) for all $x \in \mathbb{R}^n$.
For $a \in~\mathbb{R}^n, b \in~\mathbb{R}_+$, the ball $\mathcal{B}(a, b)$ is defined as  $\mathcal{B}(a, b) = \{x : \norm{a - x} \leq b \}$.
Finally, for a state $x \in \mathbb{R}^n$, we let $x^*$ denote its equilibrium value.

\section{Problem formulation}\label{sec:Problem_formulation}

\subsection{Network model}\label{sec:Network_model}

We describe the power network by a connected graph $(\mathcal{N},\mathcal{E})$ where $\mathcal{N} = \{1,2,\dots,|\mathcal{N}|\}$ is the set of buses and $\mathcal{E} \subseteq \mathcal{N} \times \mathcal{N}$ the set of transmission lines connecting the buses.
Furthermore, we use $(k,l)$ to denote the link connecting buses $k$ and $l$ and assume that the graph $(\mathcal{N},\mathcal{E})$ is directed with an arbitrary orientation, so that if $(k,l) \in \mathcal{E}$ then $(l,k) \notin \mathcal{E}$. 
For each $j \in \mathcal{N}$, we define the sets of predecessor and successor buses by  $\mathcal{N}^p_{j} = \{k : (k,j) \in \mathcal{E}\}$ and $\mathcal{N}^s_{j} = \{k : (j,k) \in \mathcal{E}\}$ respectively.
The structure of the network can be represented by its incidence matrix $H \in \mathbb{R}^{|\mathcal{N}| \times |\mathcal{E}|}$, defined as
\begin{equation*}
H_{kq} = \begin{cases}
+1, \text{ if } k \text{ is the positive end of edge } q, \\
-1, \text{ if } k \text{ is the negative end of edge } q, \\
0, \text{ otherwise.}
\end{cases}
\end{equation*}
 It should be noted that any change in the graph ordering does not alter the form of the considered dynamics. In addition, all the results presented in this paper are independent of the choice of graph ordering. 
We consider the following conditions for the power network: \newline
1) Bus voltage magnitudes are $|V_j| = 1$ per unit for all $j \in \mathcal{N}$. \newline
2) Lines $(k,l) \in \mathcal{E}$ are lossless and characterized by the magnitudes of their susceptances $B_{kl} = B_{lk} > 0$. \newline
3) Reactive power flows do not affect bus voltage phase angles and frequencies.
\newline
These conditions have been widely used in the literature in studies associated with frequency regulation, e.g. \cite{hidalgo2018frequency}, \cite{kasis2016primary}, \cite{kasis2017stability}, \cite{zhao2014design}. 
In practice, they are valid in medium to high voltage transmission systems since transmission lines are dominantly inductive and voltage variations are small.
It should be noted that all results presented in this paper are  verified with numerical simulations, presented in Section \ref{Simulation_NPCC}, on a more detailed model than our analytical one which includes voltage dynamics, line resistances and reactive power flows.

We use the swing equations to describe the rate of change of frequency at each bus. This motivates the following system dynamics (e.g. \cite{machowski2011power}),
\begin{subequations} \label{sys1}
\begin{equation}
\dot{\eta}_{kl} = \omega_k - \omega_l, \; (k,l) \in \mathcal{E}, \label{sys1a}
\end{equation}
\begin{equation}
 M^0_{j} \dot{\omega}_j \hspace{-0.25mm}=\hspace{-0.25mm} - p_j ^L + p_j^M - d^c_j - d^u_j + p^v_j - \hspace{-0.75mm}\sum_{k \in \mathcal{N}^s_j}\hspace{-0.5mm} p_{jk} +\hspace{-0.75mm} \sum_{l \in \mathcal{N}^p_j}\hspace{-0.75mm} p_{lj}, j\in \mathcal{N}, \label{sys1b}
 \end{equation}
 \begin{equation}
p_{kl}=B_{kl} \sin \eta _{kl}, \; (k,l) \in \mathcal{E}. \label{sys1c}
\end{equation}
\end{subequations}

In system~\eqref{sys1},  variables  $p^M_j$ and $\omega_j$ represent, respectively,
the mechanical power injection and
the deviation from the nominal value\footnote{We define the nominal value as an equilibrium of \eqref{sys1} with frequency equal to 50 Hz (or 60 Hz).} of the frequency at bus $j$.
Variables $d^c_j$ and $d^u_j$   represent the controllable demand and the uncontrollable frequency-dependent load and generation damping present at bus $j$ respectively.
Furthermore,
 variables $\eta_{kl}$ and $p_{kl}$ represent, respectively, the power angle difference
and the power transmitted from bus $k$ to bus $l$.
The positive constant $M^0_j$ denotes the physical inertia at bus $j$.
 Moreover, the constant $p^L_j$ denotes the frequency-independent load  at bus $j$.
Finally, the variable $p^v_j$ denotes the power injection  associated with time-varying inertia at bus $j$. 
Its dynamics follow the  virtual inertia schemes presented in e.g. \cite{lopes2014self}, \cite{d2013virtual}, \cite{wu2018frequency},
\begin{equation}\label{virtual_inertia_dynamics}
p^v_j = -M^v_j \dot{\omega}_j - D^v_j \omega_j, j \in \mathcal{N},
\end{equation}
but with a time-varying value of the virtual inertia $M^v_j$.
In particular, in \eqref{virtual_inertia_dynamics}, 
$M^v_j$ is a  non-negative  time-dependent variable describing the time-varying virtual inertia at bus $j$ and the constant $D^v_j \geq 0$ corresponds to the frequency damping coefficient associated with $p^v_j$.
\begin{remark}
We have opted to consider a constant rather than a time-varying damping coefficient $D^v_j$  in \eqref{virtual_inertia_dynamics}. 
This choice is made for simplicity and to keep the focus of the paper on the time-varying inertia. 
Including time-varying damping coefficients would result in non-existence of equilibria in primary frequency control since  these characterize the equilibrium frequency.
Jointly considering both time-varying inertia and time-varying damping coefficients is an interesting research problem for future work.
\end{remark}


\an{
\begin{remark}\label{remark_model}
The study of controllable inertia in power systems is also relevant when considering faster dynamics, such as electromagnetic models.
We study the interaction of time-varying inertia with the power system on a swing equation model for two reasons.
First, to demonstrate that the interaction of time-varying inertia with the frequency dynamics, and not any other mechanism, may yield unstable behaviour (see Section \ref{Sec_Instability}).
Second,  passivity theory, a key aspect of our analysis, is applicable to various power system models, including electromagnetic models (see e.g. \cite{trip2016internal}). Hence, it can be shown that the presented results are applicable to a range of suitable models.
\end{remark}
}

It will be convenient to define the time-varying parameters $M_j = M^0_j + M^v_j$ describing the aggregate inertia at bus $j$.
In addition, we consider the net supply variables $s_j$, defined as the aggregation of the mechanical power supply, the controllable demand, the uncontrollable frequency-dependent load and the generation and virtual inertia damping present at bus $j$, as given below
\begin{equation} \label{ssys} 
s_j = p^M_j - d^c_j  - d^u_j - D^v_j \omega_j, \; j \in \mathcal{N}.
\end{equation}

The above enable to compactly represent  \eqref{sys1}-\eqref{ssys} by
\begin{subequations} \label{sys1_compact}
\begin{align}
\dot{\eta} &= H^T \omega, \label{sys1a_compact}
\\
 M \dot{\omega}&= - p^L + s - Hp, \label{sys1b_compact}
\\
p &= B   \text{ \textbf{sin}}  (\eta), \label{sys1c_compact}
\end{align}
\end{subequations}
where  $\eta, p \in \mathbb{R}^{|\mathcal{E}|}$ and $\omega, s, p^L \in \mathbb{R}^{|\mathcal{N}|}$ are vectors associated with variables $\eta_{kl}, p_{kl}, (k,l) \in \mathcal{E}$ and $\omega_j, s_j, p^L_j, j \in \mathcal{N}$ respectively. Furthermore, $M \in \mathbb{R}^{|\mathcal{N}| \times |\mathcal{N}|}$ and $B \in \mathbb{R}^{|\mathcal{E}| \times |\mathcal{E}|}$ are diagonal matrices containing  the variables $M_j, j \in \mathcal{N}$ and parameters $B_{kl}, (k,l) \in \mathcal{E}$.

\subsection{Power supply dynamics}\label{Sec: Generation and demand dynamics}

\ak{Power supply variables represent the aggregation of the mechanical generation, controllable demand, and uncontrollable frequency dependent demand and frequency damping, as follows from \eqref{ssys}.}
To investigate a broad class of power supply dynamics,  
we will consider the following general dynamic description 
\begin{equation} \label{dynsys}
\begin{aligned}
&\dot{x}^s_j = f_j(x^s_j,-\omega_j),\\
&s_j = g_j(x^s_j,-\omega_j),
\end{aligned} \hspace{2.5mm} j \in \mathcal{N},
\end{equation}
 where $x^s_j \in \mathbb{R}^{n_j}$ denotes the internal states of the power supply variables used to update the outputs $s_j, j \in \mathcal{N}$.
 In addition, we assume that the maps $f_j : \mathbb{R}^{n_j} \times \mathbb{R} \to \mathbb{R}^{n_j}$ and $g_j : \mathbb{R}^{n_j} \times \mathbb{R} \to \mathbb{R}$ for all $j \in \mathcal{N}$ are locally Lipschitz continuous.
Moreover, we assume that in~\eqref{dynsys},
for any constant input $\omega_j(t) = \bar{\omega}_j$, there exists a unique
locally asymptotically stable equilibrium point $\bar{x}^s_j \in \mathbb{R}^{n_j}$, i.e. satisfying $f_j(\bar{x}^s_j, -\bar{\omega}_j) = 0$.
The region of attraction of~$\bar{x}^s_j$ is denoted by $\Psi_j$. 
Moreover, for any input $\omega_j \rightarrow \bar{\omega}_j$, any solution $x^s_j(t)$ of~\eqref{dynsys} such that $x^s_j(t)\in \Psi_j, t\ge 0$ must satisfy $x^s_j(t) \to \bar{x}^s_j$ as $t \to \infty$.
To facilitate the characterization of the equilibria, we also define the static input-state characteristic  map $k_{x,j} : \mathbb{R} \to \mathbb{R}^{n_j}$ as
$
k_{x,j}(-\bar{\omega}_j) := \bar{x}^s_j, j \in \mathcal{N},
$
such that $f_j(k_{x,j}(-\bar{\omega}_j), -\bar{\omega}_j) = 0$.

It should be noted that the dynamics in \eqref{dynsys} are decentralized, depending only on the local frequency $\omega_j$ for each $j \in \mathcal{N}$.
For notational convenience, we collect the variables in~\eqref{dynsys} into the vector $x^s = [x{^{s}_j}]_{j \in \mathcal{N}}$.

\begin{remark}
\label{remark_power_supply}
\ak{The quantities comprising the power supply dynamics (mechanical generation, controllable demand, uncontrollable frequency dependent demand and frequency damping) include their own dynamics} and could be represented in analogy to \eqref{dynsys}. We opted to consider a combined representation of these quantities for simplicity in presentation. However, the results presented in the paper can be trivially extended to the case where these quantities are described as individual dynamical systems.
\ak{It should also be noted that the dynamics in \eqref{dynsys} are defined with $-\omega_j$ as their input, rather than just $\omega_j$, to enable the use of passivity-based conditions (see Section \ref{Passive Dynamics}).
Such conditions have been extensively used in the literature, on network related studies (see e.g. \cite{kasis2016primary,wen2004unifying}).}
\end{remark}


\subsection{Problem statement}\label{sec_Problem}

\an{In Section \ref{Simulation_NPCC}, we present numerical simulations that  demonstrate how varying inertia may induce large frequency oscillations and  compromise the stability of the power network.
Hence,}
this study aims to  provide  local analytic conditions that associate the inertia variations and power supply dynamics such that stability is guaranteed.
The problem is stated as follows.

\begin{problem}\label{problem_definition}
Provide conditions on the time-varying inertia and power supply dynamics associated with \eqref{sys1_compact}--\eqref{dynsys} that:
\begin{enumerate}[(i)]
\item Enable asymptotic stability guarantees. 
\item Are locally verifiable.
\item Apply to high order linear and nonlinear power supply dynamics.
\item Are independent of the (connected) network topology.
\end{enumerate}
\end{problem}

The first aim requires conditions that enable asymptotic stability guarantees for the power system. 
The second objective requires conditions that can be verified using local information,  enabling plug-and-play designs. 
In addition, to enhance the practicality of our results, it is desired that those include a broad range of power supply  dynamics, including high order linear and nonlinear dynamics.
Lastly, we aim for conditions that are applicable to general network topologies, i.e.  that do not rely on the power network structure, to enable scalable designs.

\an{
In addition to Problem \ref{problem_definition}, we aim to analytically examine how time-varying inertia may lead to instability.
In particular, we intend to determine how the violation of the stability conditions established in solving Problem \ref{problem_definition} may lead to unstable behaviour in power systems.
This  analysis will theoretically validate the  inertia-induced instability phenomena observed through simulations in Section \ref{Simulation_NPCC}, and enhance our understanding on the mechanisms that cause unstable behaviour.
}

\section{Conditions on power supply  dynamics}\label{sec:conditions}

In this section we study the equilibria of \eqref{sys1_compact}--\eqref{dynsys} and provide analytic conditions on the power supply dynamics  that are subsequently used to solve Problem \ref{problem_definition}.

\subsection{Equilibrium analysis}

We now define the equilibria of the system~\eqref{sys1_compact}--\eqref{dynsys}.

\begin{definition} \label{eqbrdef}
The constants $(\eta^*, \omega^*, x^{s,*})$ define an equilibrium of the system~\eqref{sys1_compact}--\eqref{dynsys} if the following hold
\begin{subequations} \label{eqbr}
\begin{align}
0 &= H^T\omega^*, \label{eqbr1}
\\
0 &= - p^L  + s^* - Hp^*, \label{eqbr2} \\
x^{s,*}_j &= k_{x,j} (-\omega^*_j), \; j \in \mathcal{N}, \label{eqbr3} 
\end{align}
where  $p^*$ and $s^*$ in~\eqref{eqbr2} are given by
\begin{align}
\hspace{7mm} p^* &=B \text{ \textbf{sin}} (\eta^*), \label{eqbr4} \\
s^{*}_j &= g_j(x^{s,*}_j, -\omega^*_j), \; j \in \mathcal{N}. \label{eqbr5} 
\end{align}
\end{subequations}
\end{definition}

For compactness in presentation, we let $\beta = (\eta, \omega, x^s)$ where $\beta \in \mathbb{R}^{m}, m = |\mathcal{E}| + |\mathcal{N}| + \sum_{j \in \mathcal{N}} n_j$. 

\begin{remark}\label{nonunique_eq}
The equilibrium frequency $\omega^*$ uniquely defines the values of $x^{s,*}$ and $s^*$ due to the uniqueness property of the static input-state  maps $k_{x,j}, j \in \mathcal{N}$ described in Section~\ref{Sec: Generation and demand dynamics}.
By contrast, the equilibrium values of $\eta^*$ and correspondingly $p^*$ are not, in general, unique. However, these are unique under specific network configurations, such as in tree networks.
\end{remark}

\begin{remark}\label{remark_inertia}
It should be noted that the time-dependent inertia $M_j, j \in \mathcal{N},$ does not appear in the equilibrium conditions. The latter follows directly from \eqref{sys1b_compact}, i.e. the inertia affects the rate of change of frequency but not its equilibrium value.
However, as shall be discussed in the following sections, the inertia trajectories have significant impact on whether the system will converge to an equilibrium; i.e., they determine the stability properties of the equilibria.
\end{remark}

It should be noted that it trivially follows from \eqref{eqbr1} that the equilibrium frequencies of \an{\eqref{sys1_compact}, \eqref{dynsys}} synchronize, i.e. they satisfy $\omega^*_i = \omega^*_j = \omega^{s,\ast}, \forall i, j \in \mathcal{N}$, where $\omega^{s,\ast} \in \mathbb{R}$ denotes their common value.

For the remainder of the paper we assume the existence of some equilibrium to \eqref{sys1_compact}--\eqref{dynsys} following Definition \ref{eqbrdef}.
Any such equilibrium is described by $\beta^* = (\eta^*, \omega^*, x^{s,*})$. 
Furthermore, we make the following assumption
on the equilibrium power angle differences.
\begin{assumption} \label{assum1}
$| \eta^*_{ij} | < \tfrac{\pi}{2}$ for all $(i,j) \in \mathcal{E}$.
\end{assumption}

The condition imposed by Assumption \ref{assum1} can be interpreted as a security constraint that enables to deduce local convergence. In addition, it is associated with the existence of a synchronizing frequency  (see~\cite{dorfler2013synchronization}).




\subsection{Passivity conditions on power supply dynamics} \label{Passive Dynamics}

In this section we impose conditions on the power supply dynamics which  will be used to prove our main convergence result in Section \ref{Sec: Convergence}. In particular, we introduce the following passivity notion for dynamics described by \eqref{dynsys}.

\begin{definition}\label{Passivity_Definition}
System~\eqref{dynsys} is said to be locally input strictly passive with strictness constant $\rho_j$ about the constant input value $-\bar{\omega}_j$ and the constant state values~$\bar{x}^s_j$, if there exist open neighbourhoods $\Omega_j$ of $\bar{\omega}_j$ and $X_j$ of $\bar{x}^s_j$ and a continuously differentiable, positive semidefinite function $V_j(x^s_j)$ (the storage function), with a strict local minimum at $x^s_j = \bar{x}^s_j$, such that for all $\omega_j \in \Omega_j$ and all $x^s_j \in X_j$, 
\begin{equation*}
\dot{V}_j \leq (-\omega_j - (-\bar{\omega}_j)) (s_j - \bar{s}_j) - \rho_j(-\omega_j - (-\bar{\omega}_j))^2,
\end{equation*} 
where $\rho_j > 0$ 
 and $\bar{s}_j = g_j(k_{x,j}(-\bar{\omega}_j), -\bar{\omega}_j)$. 
\end{definition}

Definition \ref{Passivity_Definition} introduces an adapted notion of passivity that is suitable for the subsequent analysis\footnote{Definitions for several notions of passivity are available in \cite[Ch. 6]{khalil1996nonlinear}.}.
Passivity is a tool that has been extensively used in the literature to deduce network stability, see e.g.  \cite{kasis2016primary}, \cite{wen2004unifying}, \cite{kasis2019secondary}.
This property is easily verifiable for a wide range of systems. In particular, for linear systems it can be verified using the KYP Lemma \cite{khalil1996nonlinear} by means of a linear matrix inequality (LMI), which allows to form a convex optimization problem that can be efficiently solved. An additional approach to verify the passivity property for linear systems is to test that the corresponding Laplace transfer functions are positive real. For stable linear systems, positive realness is equivalent to the frequency response lying on the right half complex plane. 
These concepts extend  to the case of passivity with given strictness constant. 
To \ak{demonstrate the applicability of these conditions, in Section \ref{Sec: Discussion} we provide  examples of linear and nonlinear power supply dynamics that satisfy the properties presented in Definition \ref{Passivity_Definition}.} In addition, we form a suitable optimization problem that allows to deduce the storage function and the corresponding strictness constant \ak{for linear systems}.
\ak{It should further be noted that the class of dynamics that satisfy the considered passivity notion  may include fast frequency response schemes and schemes with deadband dynamics (see also \cite{kasis2017primary}).}

Below, we assume that the power supply dynamics  at each bus are locally input strictly passive with some strictness constant $\rho_j$, following Definition \ref{Passivity_Definition}.
This is a decentralized condition and hence locally verifiable, involving only the local power supply dynamics at each bus.

\begin{assumption} \label{assum2}
Each of the systems defined in~\eqref{dynsys} with input $-\omega_j$ and output $s_j$  are locally input strictly passive with strictness constant $\rho_j$, about their equilibrium values $-\omega^*_j$ and $x^{s,*}_j$  in the sense described in Definition~\ref{Passivity_Definition}.
\end{assumption}

Assumption \ref{assum2} is a key condition that allows to deduce the stability of power networks with constant inertia, as shown in \cite{kasis2016primary}.
This condition is satisfied by a wide class of power supply dynamics, including high order linear and nonlinear dynamics. Several examples of dynamics that satisfy the proposed condition are presented in Section \ref{Sec: Discussion}.
Note that, since the power supply dynamics comprise of the aggregation of the generation, controllable demand and uncontrollable frequency-dependent demand dynamics and frequency damping, Assumption \ref{assum2} allows the inclusion of dynamics that are not individually passive.

\section{Stability analysis}\label{Sec: Convergence}

In this section we provide analytic conditions on the   inertia trajectories that allow us to solve Problem \ref{problem_definition}. In addition, we provide our main stability result.
 Note that in the analysis below we study  \eqref{sys1_compact}--\eqref{dynsys} as a 
 system with states  $(\eta, \omega, x^{s})$ and time-dependent parameters $M_j(t), j \in \mathcal{N}$.

\subsection{Conditions on time-varying inertia}\label{sec_conditions_inertia}

Varying inertia may compromise the stability of the power network, as demonstrated analytically in Section \ref{Sec_Instability} and  with  simulations in Section \ref{Simulation_NPCC}.
In this section, we present conditions on the inertia trajectories that enable the study of solutions to \eqref{sys1_compact}--\eqref{dynsys} and,
as demonstrated in the following section, the provision of analytic stability guarantees.
The first condition on the inertia trajectories is presented below.

\begin{assumption}\label{assum_Virtual_Inertia}
The inertia trajectories $M_j(t), j \in \mathcal{N}$  are \ka{globally} Lipschitz in $t$ \ka{and bounded} for all $t \geq 0$.
\end{assumption}

Assumption \ref{assum_Virtual_Inertia}  requires the Lipschitz continuity \ka{and boundedness} of the inertia time-trajectories. This is a technical condition that enables the study and analysis of the solutions to \eqref{sys1_compact}--\eqref{dynsys}.
Assumption \ref{assum_Virtual_Inertia} allows to deduce the existence and uniqueness of solutions to \eqref{sys1_compact}--\eqref{dynsys}, as demonstrated in Lemma \ref{lemma_existence_uniqueness} below, proven in the Appendix.

\begin{lemma}\label{lemma_existence_uniqueness}
For any trajectory $M_j(t), j \in \mathcal{N}, t \geq 0$  that satisfies Assumption \ref{assum_Virtual_Inertia} and any initial condition $\beta(0) \in \mathbb{R}^{m}$, there exists a unique solution $\beta(t), t \geq 0$ to \eqref{sys1_compact}--\eqref{dynsys}.
\end{lemma}

The following assumption restricts the rate at which inertia trajectories may grow.

\begin{assumption}\label{assum_inertia_growth}
The inertia trajectories $M_j(t), j \in \mathcal{N}$, satisfy  $\dot{M}_j(t) \ank{\leq} 2\rho_j \ank{- \epsilon_j}, j \in \mathcal{N}$ for all $t \geq 0$, 
where $\rho_j$ is the strictness constant associated with bus $j$, in the sense described in Definition \ref{Passivity_Definition}, \ank{and $\epsilon_j, j \in \mathcal{N}$, are positive constants.}
\end{assumption}

Assumption \ref{assum_inertia_growth} restricts the rate of growth of the  inertia trajectories to be less than twice the local strictness constant $\rho_j$ associated with Assumption \ref{assum2}.
Hence, the condition relates the power supply dynamics at each bus with the rate at which inertia is allowed to grow.
Assumption \ref{assum_inertia_growth} provides a guideline for local control designs on the virtual inertia variations and could also be used from the operator's side as a means to avoid inertia induced instability.


\ak{
\begin{remark}\label{remark_inertia_growth}
Assumption \ref{assum_inertia_growth} restricts the rate of growth of inertia but not the rate at which inertia may be removed from the network, noting that inertia trajectories are lower bounded at all times by a positive constant.
This contradicts conventional intuition in power systems which suggests that high inertia values lead to enhanced stability properties.
The latter follows since higher inertia values reduce the frequency overshoot when  fault events and large disturbances occur, making it less likely that security mechanisms, which disconnect generation units and may possibly lead to blackout events, become active.
The intuition behind Assumption \ref{assum_inertia_growth} follows by noting that it bounds the rate of growth of "energy" in a Lyapunov sense, which enables to deduce stability, as demonstrated below.
\end{remark}
}

%



\subsection{Stability theorem}\label{sec_Stability_theorem}

In this section we present our main stability result  concerning the system  \eqref{sys1_compact}--\eqref{dynsys}.
In particular, the following theorem, proven in the Appendix, shows the local asymptotic convergence of solutions to \eqref{sys1_compact}--\eqref{dynsys}.

\begin{theorem}\label{thm_conv}
Let Assumptions \ref{assum_Virtual_Inertia} and \ref{assum_inertia_growth} hold and consider an equilibrium of \eqref{sys1_compact}--\eqref{dynsys} where Assumptions \ref{assum1} and \ref{assum2} hold.
Then,
there exists an open neighbourhood $\Xi$ containing that equilibrium such that solutions $\beta(t), t \geq 0$ to \eqref{sys1_compact}--\eqref{dynsys} \ank{initiated in $\Xi$} asymptotically converge to the set of equilibria \ank{of \eqref{sys1_compact}--\eqref{dynsys}}. 
\end{theorem}

Theorem \ref{thm_conv} provides analytic guarantees for the stability of the power network at the presence of time-varying  inertia.
Note that the main conditions on the power supply dynamics and the varying inertia trajectories, described by Assumptions \ref{assum2} and \ref{assum_Virtual_Inertia}--\ref{assum_inertia_growth} respectively, are locally verifiable.
These conditions may be used for the  design of prototypes and guidelines for virtual inertia  schemes and may motivate  practical control designs that enhance the stability properties of the power grid.
Furthermore, Assumption \ref{assum2} includes a wide range of power supply dynamics, as demonstrated in the following section.
In addition, Theorem \ref{thm_conv} applies to any connected power network topology.
Hence, all objectives of Problem \ref{problem_definition} are satisfied.


\ak{
\begin{remark}\label{remark_discontinous}
The convergence analysis presented above assumes that inertia trajectories are  Lipschitz continuous, through Assumption \ref{assum_Virtual_Inertia}.
An interesting and relevant extension to this work is to  consider jumps on inertia trajectories.
However, this direction raises significant challenges in terms of determining suitable bounds on jump magnitudes and frequency of occurrence. 
In addition, it is important that conditions on inertia trajectories are locally verifiable, to enable scalable, plug-and-play designs, that are of particular importance in large-scale power systems. 
Hence, this direction deserves independent research attention and is left as future work. 
\end{remark}
}

\an{
\begin{remark}\label{rem_comparison}
Existing works, such as \cite{kasis2016primary}, rely on passivity arguments applied to continuous power supply dynamics, by requiring  marginal strictness to guarantee convergence. 
In contrast, \ank{our analysis presents a systematic approach that utilizes structured deviations from passivity within local power supply dynamics, through the term $\rho_j(-\omega_j - (-\bar{\omega}_j))^2$ in Definition \ref{Passivity_Definition}, to develop design conditions for local inertia schemes, as described by Assumption \ref{assum_inertia_growth}.} 
This also enables the co-design of power supply and local inertia  dynamics.
\end{remark}
}

\section{Application examples}\label{Sec: Discussion}

In this section we provide two examples \ak{of linear and one example of nonlinear} power supply dynamics that fit within the presented framework. In addition, we present a systematic approach that allows us to obtain the strictness constant $\rho_j$ associated with Assumption \ref{assum2} for linear systems, which provides the local bound on the maximum rate of inertia growth presented in Assumption \ref{assum_inertia_growth}.

In particular, we consider general linear power supply dynamics with minimal state space realization of the form
\begin{equation}\label{sys_linear}
\begin{aligned} 
\dot{x}^s_j &= A_j x^s_j + B_j (-\omega_j), \\
s_j &= C_j x^s_j + D_j (-\omega_j),
\end{aligned}
\end{equation}
where    $A_j \in \mathbb{R}^{n_j \times n_j}, B_j \in \mathbb{R}^{n_j}, C_j \in \mathbb{R}^{1 \times n_j}$ and $D_j \in \mathbb{R}$ are matrices describing the power supply dynamics at bus $j$.
For the dynamics described in \eqref{sys_linear} it can be deduced, by suitably adapting \cite[Thm. 3]{willems1972dissipative}, that the strictness constant $\rho_j$ may be efficiently obtained as the solution of the following optimization problem:
\begin{equation}\label{LMI_Optimization_Problem}
\begin{aligned}
&\max_{\hat{D}, P} D_j -  \hat{D} \\
s.t. \quad & P = P^T \succeq 0, \\
&\begin{bmatrix}
A_j^TP + A_jP \quad PB_j - C_j^T \\
B_j^TP - C_j \quad -2\hat{D}
\end{bmatrix} \preceq 0,
\end{aligned}
\end{equation}
i.e., when \eqref{LMI_Optimization_Problem} is maximized at $\hat{D} = \hat{D}^*$ and some $P$, then $\rho_j = D_j - \hat{D}^*$.
The above problem can be solved in a computationally efficient manner, using standard semidefinite programming tools. 
In addition, the matrix $P$ can be used to obtain the storage function associated with $s_j$, as $V(x^s_j) = \frac{1}{2} (x^s_j)^TPx^s_j$. 
Furthermore, it is intuitive to note that   a larger value of $D_j$, which describes the local damping, yields a larger strictness constant $\rho_j$, as follows from \eqref{LMI_Optimization_Problem}.
Below, we present two examples of power supply dynamics and explain how the strictness constant may be obtained in each case.

As a first example, we consider the first order generation dynamics considered e.g. in  \cite{trip2017distributed}, that describe the time lag between changes in frequency and the response from generation units.
The  dynamics  are given by
\begin{equation}\label{First_order_system}
\begin{aligned}
\tau_j \dot{x}_j &= - x_j - K_j \omega_j,\\
s_j &= x_j - \lambda_j \omega_j,
\end{aligned}
\end{equation}
where $\tau_j > 0, K_j > 0$ and $\lambda_j > 0$ are the time, droop and damping constants respectively.
The solution to \eqref{LMI_Optimization_Problem} for system \eqref{First_order_system} is given by $\hat{D} = 0$ for any $\tau_j, K_j$, which results in 
$\rho_j = \lambda_j$. The latter suggests by Assumption \ref{assum_inertia_growth} that $\dot{M}_j < 2\lambda_j$ should hold.

A more involved example that demonstrates the applicability of the presented analysis concerns the fifth order turbine governor dynamics provided by the Power System Toolbox \cite{cheung2009power}.
This model is described in the Laplace domain by the following transfer function
\begin{equation}\label{fifth_order}
G_j(s)=K_j\frac{1}{(1+sT_{s,j})}\frac{(1+sT_{3,j})}{(1+sT_{c,j})}\frac{(1+sT_{4,j})}{(1+sT_{5,j})} + \lambda_j,
\end{equation}
relating the power supply output $\hat{s}_j$ with the frequency deviation input $-\hat{\omega}_j$, where $K_j$ and $T_{s,j}, T_{3,j}, T_{c,j}, T_{4,j},$ $T_{5,j}$ are the droop coefficient and time-constants respectively and $\lambda_j$ denotes the frequency damping.
Realistic values for the coefficients in \eqref{fifth_order} are provided in \cite{cheung2009power}.
To provide a numerical example on how the strictness constant $\rho_j$ may be obtained, we consider the turbine governor dynamics at bus $36$ of the NPCC network, where the above coefficients take the values $(K_j, T_{s,j}, T_{3,j}, T_{c,j}, T_{4,j}, T_{5,j}, \lambda_j) = (110.1, 0.45, 0.1, 0, 13.25, 54, 30.3)$.
By solving \eqref{LMI_Optimization_Problem} using the CVX toolbox \cite{grant2009cvx}, we obtain a strictness coefficient $\rho_j$ of approximately $28.0$.

\begin{figure}
\includegraphics[trim = 7mm 0 0 3mm, clip, scale=0.64]{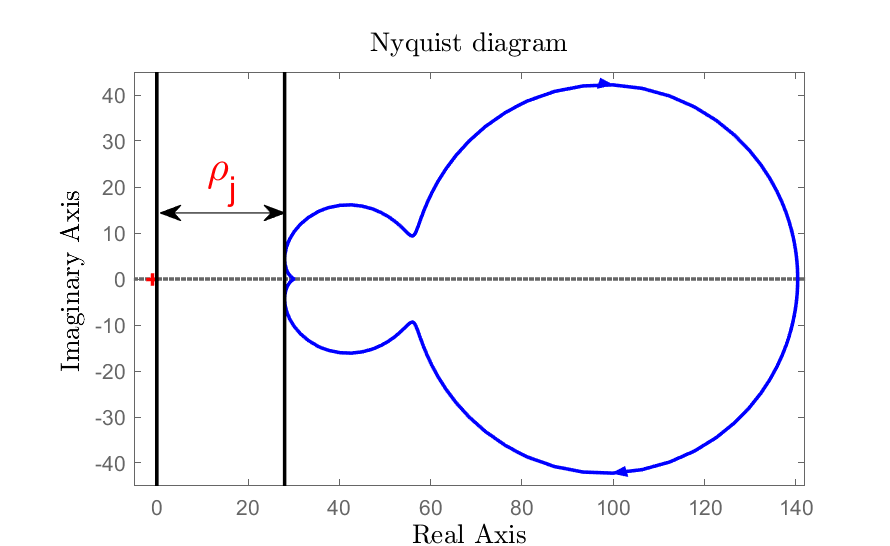}
\vspace{-2mm}
\caption{Nyquist plot for \eqref{fifth_order}, with coefficients associated with the turbine governor dynamics at bus 36 within the NPCC network. The strictness constant $\rho_j$ corresponds  to the horizontal distance between the Nyquist plot and the imaginary axis.}
\label{Fig_Nyquist}
\end{figure}

A graphical approach can also be used to obtain the strictness constant $\rho_j$ for systems described by \eqref{sys_linear}.
In particular,  an approach to verify passivity\footnote{Passive systems satisfy Definition \ref{Passivity_Definition} with strictness constant $\rho_j = 0$.} for  stable linear systems, is to test that their transfer function is positive real, which is equivalent to the frequency response lying on the right half complex plane. 
The strictness constant $\rho_j$ can be obtained as the horizontal distance between the Nyquist plot and the imaginary axis.
This is demonstrated in Fig. \ref{Fig_Nyquist}, which depicts the Nyquist plot for \eqref{fifth_order} with the coefficients given as above. The distance between the plot and the imaginary axis matches exactly the value obtained by solving \eqref{LMI_Optimization_Problem}.

\ak{It should be noted that Assumption \ref{assum2} is also applicable to nonlinear power supply dynamics. As an example, consider local power supply dynamics described by a non-decreasing function $\hat{f}_j$ and some frequency damping $\lambda_j$ as follows
\begin{equation}\label{example_static_nonlinearity}
s_j = \hat{f}_j(-\omega_j) - \lambda_j\omega_j.
\end{equation}
For system \eqref{example_static_nonlinearity}, it follows that $\rho_j = \lambda_j$. Additional results can be obtained for general nonlinear power supply dynamics, following the description in \eqref{dynsys}, using $\mathcal{L}_2$-gain arguments similar to \cite[Prop. 1]{kasis2016primary}.
}

\section{Inertia induced instability}\label{Sec_Instability}


\an{A key aspect highlighted in this paper, and demonstrated through simulations in Section \ref{Simulation_NPCC}, is that varying inertia may cause unstable behaviour. This 
motivated the analytical study of the stability properties of power systems with time-varying inertia in Section \ref{Sec: Convergence}.
 This section aims to provide a comprehensive understanding of the  inertia-induced instability issue by offering theoretical validation, generalizability, and deeper insight into the mechanisms causing instability.
 Moreover, it aims to establish that such behaviour originates from the fundamental frequency dynamics.
In particular, we provide sufficient conditions that allow us to deduce, from  any arbitrarily small deviation from the equilibrium frequency, the existence of  single-bus inertia variations that cause substantial deviations in  system trajectories.}
The latter describes unstable behaviour, see e.g. \cite{haddad2011nonlinear}.
Hence, we aim to demonstrate how systematic unstable behaviour may result due to inertia variations.

For simplicity, we shall consider linearized power flow equations for the dynamics in \eqref{sys1_compact}.
The dynamics of such a system are:
\begin{subequations} \label{sysl_compact}
\begin{align}
\dot{\eta} &= H^T \omega, \label{sysla_compact}
\\
 M \dot{\omega}&= - p^L + s - Hp, \label{syslb_compact}
\\
p &= B    \eta. \label{syslc_compact}
\end{align}
\end{subequations}
The equilibria of \eqref{dynsys}, \eqref{sysl_compact} are described in analogy to Definition \ref{eqbrdef} as follows.

\begin{definition} \label{eqbr_def_lin}
The constants $(\eta^*, \omega^*, x^{s,*})$ define an equilibrium of the system~\eqref{dynsys}, \eqref{sysl_compact} if \eqref{eqbr1}--\eqref{eqbr3} hold, 
where  in~\eqref{eqbr2}   $s^*$  is given by \eqref{eqbr5} and 
 $p^* =B \eta^*$. 
\end{definition}

Below, we demonstrate that the stability results presented in Section \ref{Sec: Convergence} extend to \eqref{dynsys}, \eqref{sysl_compact}. In particular, the following lemma, proven in the Appendix, demonstrates that Assumptions \ref{assum2}, \ref{assum_Virtual_Inertia} and \ref{assum_inertia_growth} suffice to deduce the local asymptotic convergence of solutions to \eqref{dynsys}, \eqref{sysl_compact}.

\begin{lemma}\label{lemma_conv}
Let Assumptions \ref{assum_Virtual_Inertia} and \ref{assum_inertia_growth} hold and consider an equilibrium of \eqref{dynsys}, \eqref{sysl_compact} where Assumption  \ref{assum2} holds.
Then, 
there exists an open neighbourhood $\Xi$ containing that equilibrium such that solutions $\beta(t), t \geq 0$ to \eqref{dynsys}, \eqref{sysl_compact} \ank{initiated in $\Xi$ asymptotically converge to the set of equilibria of \eqref{dynsys}, \eqref{sysl_compact}}. 
\end{lemma}

\ak{Lemma \ref{lemma_conv} follows similarly to Theorem \ref{thm_conv} but for the relaxed system \eqref{dynsys}, \eqref{sysl_compact}. Compared to Theorem \ref{thm_conv}, Lemma \ref{lemma_conv} does not need to impose Assumption \ref{assum1}, since linearized power transfer equations are considered, as follows from \eqref{syslc_compact}. Nevertheless, the result remains local since the constant input asymptotic stability condition on \eqref{dynsys} and Assumption \ref{assum2}  hold locally. 
}


\subsection{Conditions for instability}

To facilitate the analysis, we define the notion of $\gamma$-points. 
Such points result by considering the equilibria of \eqref{dynsys}, \eqref{sysl_compact} when some bus $k$ has a fixed frequency at all times.  
The notion of $\gamma$-points will be used to provide conditions that characterize the solutions to \eqref{dynsys}, \eqref{sysl_compact}.

\begin{definition} \label{eqbrdef_2}
The constants $(\hat{\eta}, \hat{\omega}, \hat{x}^{s})$ define a $\gamma$-point of \eqref{dynsys}, \eqref{sysl_compact} associated with fixed frequency $\bar{\omega} \in \mathbb{R}$ at some bus $k$ if the following hold
\begin{subequations} \label{eqbr_2}
\begin{align}
\hat{\omega}_k &= \bar{\omega}, \label{eqbrk_2} \\
0 &= H^T\hat{\omega}, \label{eqbr1_2}
\\
0 &=  - p^L_j + \hat{s}_j \hspace{-0.5mm}-\hspace{-1mm} \sum_{\bar{k} \in \mathcal{N}^s_j}\hspace{-0.5mm} \hat{p}_{j\bar{k}} +\hspace{-1mm} \sum_{l \in \mathcal{N}^p_j}\hspace{-0.5mm} \hat{p}_{lj}, j\in \mathcal{N}\setminus \{k\}, \label{eqbr2_2} \\
\hat{x}^{s}_j &= k_{x,j} (-\hat{\omega}_j), \; j \in \mathcal{N}, \label{eqbr3_2} 
\\ 
\hat{p} &= B \hat{\eta}, \label{eqbr4_2} \\
\hat{s}_j &= g_j(\hat{x}^{s}_j, -\hat{\omega}_j), \; j \in \mathcal{N}. \label{eqbr5_2} 
\end{align}
\end{subequations}
The set of points $(\hat{\eta}, \hat{\omega}, \hat{x}^{s})$ that satisfy \eqref{eqbr_2} associated with fixed   frequency  $\bar{\omega}$ is denoted by $\an{\Gamma_k}(\bar{\omega})$.
\end{definition}


It is intuitive to note that the notion of $\gamma$-points follows by considering the theoretical case where the inertia at a single bus $k$ within system \eqref{dynsys}, \eqref{sysl_compact} is infinite. The equilibrium points of such a system, when the initial conditions satisfy $\omega_k(0) = \bar{\omega}$, are described by  \eqref{eqbr_2}. 
Alternatively, $\gamma$-points may describe the equilibria of \eqref{dynsys}, \eqref{sysl_compact} for some $p^L$.
It should be noted that for any set $\an{\Gamma_k}(\bar{\omega})$ the values of $(\hat{\omega}, \hat{x}^s)$ are unique and satisfy $(\hat{\omega}_j, \hat{x}^s_j) = (\bar{\omega}, k_{x,j} (-\bar{\omega}_j)), j \in \mathcal{N}$.

The following assumption is the main condition imposed to deduce inertia induced instability. 

\begin{assumption}\label{assum_instability}
The following  
hold for \eqref{dynsys}, \eqref{sysl_compact}, some bus $k$ 
and some positive constants $\hat{\epsilon}, \bar{\epsilon}, \Phi$ and $\bar{\Phi}$ satisfying $\hat{\epsilon} < \bar{\epsilon} < \Phi$ and  $\bar{\Phi} = \Phi +  \bar{\epsilon}$:
\begin{enumerate}[(i)]
\item For any $\bar{\omega} \in \mathcal{B}( \omega^{s,\ast}, \bar{\Phi})$, there exist \an{positive constants  $\epsilon, {\epsilon}', \hat{\tau}$, satisfying $\epsilon < {\epsilon}' < \hat{\epsilon}$} such that when $\an{\omega_k(t) = \bar{\omega}}$, 
$\omega_j(t) \in \mathcal{B}(\bar{\omega}, \epsilon), j \in \mathcal{N}\an{\setminus \{k\}},$ 
for all $t \in [0, \hat{\tau}]$, then $\beta(\hat{\tau}) \in \mathcal{B}(\gamma, \an{{\epsilon}'}), \gamma \in \an{\Gamma_k}(\bar{\omega})$.
\item When $M^v_k(t) = 0, t \geq 0$ and 
Assumption \ref{assum_inertia_growth} holds
then, for any $\bar{\omega} \in \mathcal{B}( \omega^{s,\ast}, \bar{\Phi}) \setminus \{\omega^{s,\ast}\}$ and any solution to  \eqref{dynsys}, \eqref{sysl_compact}  there exists $\tau$ such that $\beta(0) \in \mathcal{B}(\gamma, \hat{\epsilon}) \setminus \an{\Gamma_k}(\omega^{s,\ast}), \gamma \in \an{\Gamma_k}(\bar{\omega})$ implies that $|\omega_k (\tau) - \omega^{s,\ast}| > \bar{\omega} + \bar{\epsilon}$.
\item  Assumption \ref{assum2} holds for all points in $\Upsilon = \{(\tilde{\omega}, \tilde{x}^s) : (\tilde{\omega}, \tilde{x}^s) \in \an{\Gamma_k}(\bar{\omega}), \bar{\omega} \in \mathcal{B}( \omega^{s,\ast}, \bar{\Phi})\}$.
 In addition,  for all points in $\Upsilon$, all neighbourhoods $\Omega_j$ of $\tilde{\omega}_j$ and $X_j$ of $\tilde{x}^s_j$ associated with Definition \ref{Passivity_Definition}, and regions of attraction $\Psi_j$ associated with $\tilde{x}^s_j$ in the description of \eqref{dynsys},
satisfy $\Omega_j \times X_j \supseteq \bar{\Omega}_j \times \bar{X}_j   , \Psi_j  \supseteq \bar{X}_j$, where $\bar{\Omega}_j := \{p : p \in \mathcal{B}(\omega^{s,\ast}, \bar{\Phi}) \}$
and $\bar{X}_j := \{ p : p \in \mathcal{B}(\bar{x}^s_j,  \bar{\Phi}), \bar{x}^s_j \in  \an{\Gamma_k}( \omega^{s,\ast}) \}, j \in \mathcal{N}$.
\end{enumerate}
\end{assumption}

Assumption \ref{assum_instability} is split in three parts.
Part (i) requires that when all frequencies lie in a ball of size $\epsilon$ around $\bar{\omega}$ \an{and $\omega_k = \bar{\omega}$ (a condition associated with the theoretical case where the inertia at bus $k$ is infinite),} then the solutions of the system  converge to a ball of size $\hat{\epsilon}$ around a $\gamma$-point in $\an{\Gamma_k}(\bar{\omega})$ within some finite time $\hat{\tau}$. 
This assumption is associated with Lemma \ref{lemma_conv} which states the convergence of the solutions to \eqref{dynsys}, \eqref{sysl_compact}.
Assumption \ref{assum_instability}(i) is a  mild condition that is expected to hold for almost all practical power systems.

Assumption \ref{assum_instability}(ii) is the most important condition imposed, requiring that solutions initiated at any non-equilibrium point within a ball of size $\hat{\epsilon}$ from a point within $\an{\Gamma_k}(\bar{\omega})$ will be such that the frequency deviation from equilibrium at   some given bus $k$ is larger in magnitude than $\bar{\omega} +  \bar{\epsilon}, \bar{\epsilon} > \hat{\epsilon}$ at some time $\tau$.
This condition is important as it enables the main arguments in our instability analysis.
We demonstrate with simulations in Section \ref{Simulation_NPCC} that Assumption \ref{assum_instability}(ii) applies to two realistic networks.

Finally, Assumption \ref{assum_instability}(iii) requires that the local  passivity and asymptotic stability properties on power supply dynamics associated with Assumption \ref{assum2} and the description below \eqref{dynsys} hold for a broad range of points, i.e. for all points in $\Upsilon$.
In addition, it requires sufficiently large regions where these local conditions hold.
 Assumption \ref{assum_instability}(iii) could be replaced by the simpler, but more conservative, condition that Assumption \ref{assum2} and the asymptotic stability properties on power supply dynamics hold globally for all points in $\Upsilon$.
In addition, it could further be relaxed by letting $\Upsilon = \mathbb{R}^{|\mathcal{N}| + \sum_{j \in \mathcal{N}} n_j}$.


\subsection{Instability theorem}

In this section we present our main instability results. 
In particular, the following theorem, proven in the Appendix, demonstrates the existence of single-bus inertia trajectories that cause substantial frequency deviations from any non-equilibrium initial condition at  bus $k$ frequency. The latter suggests unstable behaviour, as  shown in Corollary~\ref{cor_instability} below.

\begin{theorem}\label{thm_instability}
Let Assumption \ref{assum_Virtual_Inertia} hold, 
Assumption \ref{assum_inertia_growth} hold for all $j \in \mathcal{N} \setminus \{ k\}$
 and consider an equilibrium of \eqref{dynsys}, \eqref{sysl_compact} 
where Assumption \ref{assum_instability} holds for some bus $k$.
Then,   for 
every  $\delta > 0$ 
there exists some finite time $\tilde{t}$ and some finite trajectory  $M_k (t) \geq M^0_k, t \in [0, \tilde{t}]$ such that $|\omega_k(0) - \omega^{s,\ast}| \geq \delta$ implies   that 
 $\beta(\tilde{t}) \notin \mathcal{B}(\gamma, \Phi)$ for any  $\gamma \in \an{\Gamma_k}(\omega^{s,\ast})$
 where 
 $\Phi$ follows from Assumption \ref{assum_instability}.
\end{theorem}

Theorem \ref{thm_instability} demonstrates the existence of  single-bus inertia  trajectories such that an arbitrary small frequency deviation may result in substantial deviations in system trajectories.
 Hence, the stability properties of the power network may be compromised when local inertia is suitably controlled.
The latter may motivate attacks on the inertia of the power network, which may cause large frequency oscillations.
In addition, Theorem \ref{thm_instability} demonstrates the importance of restricting inertia trajectories at all buses, as follows from   Assumption \ref{assum_inertia_growth}, to deduce stability.
It should also be noted that a result that trivially follows from Theorem \ref{thm_instability} is the existence of inertia trajectories on multiple buses that cause instability.


The following result, proven in the Appendix, is a corollary of Theorem \ref{thm_instability} which deduces the existence of single-bus inertia trajectories that render an equilibrium point unstable.

\begin{corollary}\label{cor_instability}
Let Assumption \ref{assum_Virtual_Inertia} hold,
Assumption \ref{assum_inertia_growth} hold for all $j \in \mathcal{N} \setminus \{ k\}$
 and consider an equilibrium of \eqref{dynsys}, \eqref{sysl_compact} 
where Assumption \ref{assum_instability} holds for some bus $k$.
 Then, there exists a finite trajectory  $M_k (t) \geq M^0_k, t \geq 0$ such that the considered equilibrium is unstable.
\end{corollary}

\ak{
\begin{remark}\label{rem_instability_discussion}
The presented analysis challenges the conventional intuition that increased inertia yields enhanced  stability properties in power networks.
In addition, it showcases the importance of the results presented in Section \ref{Sec: Convergence}, concerning the stability of power networks under varying inertia.
Moreover, Assumption \ref{assum_instability} and the proof of Theorem \ref{thm_instability} offer intuition on how local inertia trajectories that lead to instability may be constructed.
\an{In particular, instability arises from the coupling between the frequency and power supply dynamics at given inertia trajectories.
If inertia increases during significant frequency deviations, then reaching the equilibrium frequency is delayed. This delay causes excess generation, leading to frequency overshoots when inertia suddenly drops, ultimately causing growing oscillations.}
The effectiveness of these results is validated in Section \ref{Simulation_NPCC}, which demonstrates single-bus instability in two  power network settings.
\end{remark}
}

\an{
\begin{remark}\label{rem_lemma_thm_comparison}
Theorem \ref{thm_instability} incorporates all  assumptions used to establish stability in Lemma \ref{lemma_conv}, except Assumption \ref{assum_inertia_growth} at bus $k$. Specifically, Theorem \ref{thm_instability} includes Assumption \ref{assum_Virtual_Inertia}, Assumption \ref{assum_inertia_growth} for all buses except bus $k$, and Assumption 5(iii), which relates to Assumption \ref{assum2}. This  comparison between the conditions for stability and instability, underscores the critical importance of Assumption~\ref{assum_inertia_growth}, providing a strong justification for its inclusion in inertia scheme design.
\end{remark}
}

\section{Simulations} \label{Simulation_NPCC}

In this section we present numerical simulations using the Northeast Power Coordinating Council (NPCC) $140$-bus system and the IEEE New York / New England $68$-bus  
system
 that further motivate and validate the main findings of  this paper.
In particular, we first demonstrate how varying virtual inertia may induce large oscillations in the power network.
 We then verify our  analytic stability results by showing that the main imposed conditions, given by Assumptions \ref{assum2} and \ref{assum_inertia_growth}, yield a stable behaviour. 
 Finally, we demonstrate how single-bus inertia variations may yield large frequency oscillations, verifying Theorem \ref{thm_instability}.

For our simulations,  we used the Power System Toolbox~\cite{cheung2009power} on Matlab.
The model used by the toolbox is more detailed than our analytic one, including  voltage dynamics, line resistances and  a transient reactance generator model\footnote{The simulation details can be found in the data files datanp48 (NPCC network) and data16em (New York / New England network) and the Power System Toolbox  manual \cite{cheung2009power}.}.

\subsection{Simulations on the NPCC network} \label{subsection_NPCC}

The NPCC network  consists of 47 generation and 93 load buses and has a total real power of $28.55$ GW. 
For our simulation, we considered a step increase in demand of magnitude $2$ p.u. (base $100$ MVA) at load buses $2$ and $3$  at $t=1$ second. The simulation precision was set at $10$ ms.

\subsubsection{Inertia induced oscillations}\label{sec_Inertia_Oscillations}

To demonstrate that varying inertia  may induce large frequency oscillations and hence compromise the stability of the power network, we considered two cases: 
\begin{enumerate}[(i)]
\item   no presence of virtual inertia,
\item  the presence of virtual inertia at $10$ generation buses (buses $23, 48, 50, 54, 56, 57, 72, 80, 82$ and $133$) of magnitude $M_a$, where $M_a$ was equal to $50\%$ of the physical inertia at bus $133$.
\end{enumerate}
The trajectories of the virtual inertia associated with case (ii)  were coupled with the frequency dynamics as follows:
\begin{equation}\label{unstable_VI}
M^v_j(t) = \begin{cases}
M_a, \text{ if } \omega_m(t) > 0.02 \text{ Hz,} \\
0, \text{ otherwise},
\end{cases}
\end{equation}
where $\omega_m(t) = \max_{j \in \mathcal{N}} |\omega_j(t)|$.
 The scheme in \eqref{unstable_VI} adds inertia to the power system when a noticeable frequency deviation is experienced and removes it when the system returns close to the nominal frequency.

\begin{figure}
\includegraphics[scale=0.62]{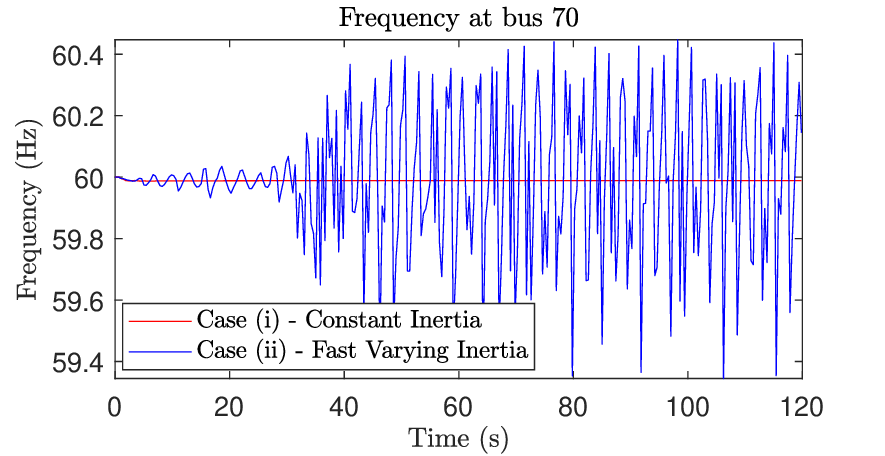}
\vspace{-3mm}
\caption{Frequency response at bus $70$ when: (i) no virtual inertia is present and (ii) virtual inertia  described by \eqref{unstable_VI} is included at $10$ buses.}
\label{Figure_Instability}
\end{figure}

The frequency response at bus $70$ for the two considered cases is presented in Fig. \ref{Figure_Instability}. From Fig. \ref{Figure_Instability}, it can be seen that the addition of fast varying inertia yields large oscillations in the power network. 
The oscillations follow due to the coupling between the frequency and generation dynamics.
In particular, generators respond to frequency signals by appropriately adapting the generated power.
However, when the inertia abruptly increases under a noticeable frequency deviation, following \eqref{unstable_VI},  it takes longer for frequency to reach its steady state.
The latter causes excess generation to be produced which induces frequency overshoots when the  inertia suddenly drops. 
The above process results in frequency oscillations, as verified in Fig. \ref{Figure_Instability}.
Therefore, care needs to be taken when varying virtual inertia is introduced in the power network,  particularly when its dynamics are coupled with those of the network.

It should be noted that the possibly destabilising effects of varying inertia are visible even when a large and robust system, such as the NPCC network, is considered. The latter highlights the potential impact of virtual inertia schemes and provides further motivation for their proper regulation.
In addition, note that we opted not to present the (widely acknowledged) benefits of (constant) virtual inertia for compactness and to keep the focus  on the impact of varying virtual inertia. 

%

\subsubsection{Stability preserving varying inertia}\label{sec_stab_varying_inertia}

To demonstrate the validity and applicability of the proposed conditions,  we repeated the above simulation with the virtual inertia satisfying Assumption \ref{assum_inertia_growth}.
In particular, we introduced varying virtual inertia of maximum magnitude $M_a$  in the same set of generation buses as in case (ii).
To comply with Assumption \ref{assum_inertia_growth}, the coupling between the virtual inertia dynamics and the frequency dynamics was given by:
\begin{subequations}\label{Stable_VI}
\begin{align}
\dot{M}^v_j &= 
\min(\tau (-M^v_j + u_j), 2\rho_j - \epsilon),\label{Stable_derivative}
\\
u_j(t) &= \begin{cases}
M_a, \text{ if } \omega_m(t) > 0.02 \text{ Hz,} \\
0, \text{ otherwise},
\end{cases} \label{set_point}
\end{align}
\end{subequations}
where  $\tau  = 100$ $\text{s}^{-1}$ was the time constant of the virtual inertia dynamics, selected such that a fast inertia variation  is allowed, and $u_j$ an input set point.
In addition,  $\rho_j$ corresponded to the strictness constant at each bus, calculated following the approach presented in Section \ref{Sec: Discussion}, and $\epsilon = 10^{-4}$ a small constant introduced to ensure that the inequality in Assumption \ref{assum_inertia_growth} was satisfied.
 The scheme in \eqref{Stable_VI} enabled fast variations in virtual inertia by setting a large value for $\tau$ and simultaneously restricted its rate of growth in accordance with Assumption \ref{assum_inertia_growth}.
 To couple the frequency and inertia dynamics, the input $u_j$ was  set to $M_a$ when the frequency deviation exceeded $0.02$~Hz and $0$ otherwise, as follows from \eqref{set_point}.
 \an{It should be noted that  $u_j$ in \eqref{set_point} has identical dynamics with $M^v_j$ in \eqref{unstable_VI}. 
This design choice, coupled with the selection of a large value for $\tau$ in \eqref{Stable_derivative}, was intended to enable \eqref{Stable_VI} to simultaneously approximate \eqref{unstable_VI} and satisfy Assumption \ref{assum_inertia_growth}.}
This case will be referred to as case (iii).

\begin{figure}[t!]
\includegraphics[scale=0.6]{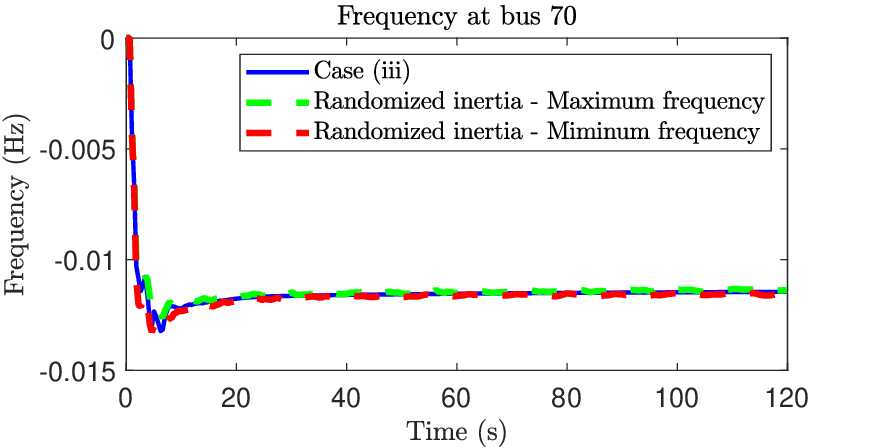}
\vspace{-5mm}
\caption{Frequency response at bus $70$ when: (i) case (iii) is implemented, with virtual inertia satisfying \eqref{Stable_VI}, (ii) the randomized scheme described by \eqref{Stable_derivative}, \eqref{randomized_set_points} is implemented, where green and red lines correspond to the maximum and minimum frequencies obtained after $500$  trials respectively.}
\label{Figure_Stability}
\vspace{-2mm}
\end{figure}

The frequency response at bus $70$ resulting from implementing \eqref{Stable_VI} is depicted in Fig. \ref{Figure_Stability}. From Fig. \ref{Figure_Stability}, it follows that the proposed scheme yields a stable response for the power system, which validates the main analytical results of the paper. 
\an{Consequently, it verifies that the resulting stable behavior arises from the introduction of the derivative bound associated with Assumption~\ref{assum_inertia_growth} in~\eqref{Stable_VI}.}

\begin{figure}[t!]
\includegraphics[trim = 3mm 0 0 0, clip, scale=0.63]{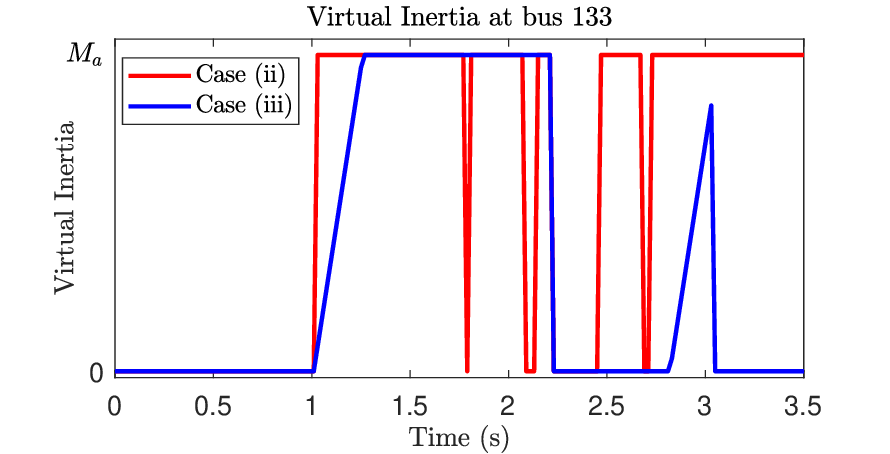}
\vspace{-5mm}
\caption{Virtual inertia at bus $133$ for cases (ii) and (iii) concerning  inertia described by \eqref{unstable_VI} and  \eqref{Stable_VI} and leading to oscillatory and stable responses respectively. After $3.5$ seconds, the virtual inertia is $M_a$ at almost all times, i.e. with fast fluctuations, for case (ii) and zero for case (iii).}
\label{Figure_VI}
\vspace{-2mm}
\end{figure}

To demonstrate that Assumption \ref{assum_inertia_growth} enables fast changes in inertia, the virtual inertia  at bus $133$ concerning cases (ii) and (iii) is depicted in Fig. \ref{Figure_VI}. 
From Fig. \ref{Figure_VI}, it follows that for case (iii),  the maximum virtual inertia support is provided within $0.3$ seconds from the time the frequency overpasses $0.02$ Hz and is removed completely after $3.5$ seconds.
By contrast, in case (ii) the virtual inertia instantly reaches $M_a$ at $1$ second but fluctuates between $0$ and $M_a$. Its fast fluctuations create frequency oscillations  which lead to more  inertia fluctuations due to \eqref{unstable_VI}, yielding the oscillatory frequency response depicted in Fig. \ref{Figure_Instability}.

To further demonstrate that the proposed conditions yield a stable response, we considered randomly changing inertia input set points described by:
\begin{equation}\label{randomized_set_points}
u_j(t^+) \hspace{-0.5mm}=\hspace{-0.5mm} \begin{cases}
u_j(t), \text{ if}~\Mod(t, 0.5) \neq 0, \\
u_j(t) \hspace{-0.5mm}+\hspace{-0.5mm} 0.5M_a, \text{ if}~\Mod(t, 0.5) = 0, r_j(t) \geq 0.5,\\
\max(u_j(t) \hspace{-0.5mm}-\hspace{-0.5mm} 0.5M_a, 0), \text{ otherwise,}
\end{cases}
\end{equation}
where $t^+ = \lim_{\epsilon \rightarrow 0} (t + \epsilon)$ and $r_j(t)$ is randomly selected  from the uniform distribution $[0,1]$ at each $t$ that satisfies $\Mod(t, 0.5) = 0$.
The scheme in \eqref{randomized_set_points} updates the set points $u_j$ every $0.5$ seconds by equiprobably increasing or decreasing their values by $0.5M_a$. Simultaneously, it ensures that the set points take non-negative values.
The dynamics \eqref{Stable_derivative}, \eqref{randomized_set_points} were implemented and simulated $500$ times on the same setting as cases (ii) and (iii), to show that stability is preserved for a wide range of varying virtual inertia profiles that satisfy Assumption \ref{assum_inertia_growth}.
The latter is demonstrated in Fig. \ref{Figure_Stability}, which shows that the maximum and minimum frequencies obtained with the presented randomized scheme are very close to the response associated with case (iii).

\ak{
To explore the conservativeness of Assumption \ref{assum_inertia_growth}, we performed simulations using the same setting described in case (iii) but with virtual inertia dynamics described by
\begin{subequations}\label{Stable_VI_Conserv}
\begin{align}
\dot{M}^v_j &= 
\min(\tau (-M^v_j + u_j), 2 \kappa\rho_j),\label{Stable_derivative_Conserv}
\\
u_j(t) &= \begin{cases}
\bar{M}_a, \text{ if } \omega_m(t) > \bar{\omega} \text{ Hz,} \\
0, \text{ otherwise},
\end{cases} \label{set_point_Conserv}
\end{align}
\end{subequations}
where $\kappa > 1$ denotes the factor by which the main condition in Assumption \ref{assum_inertia_growth} is relaxed.
In addition, $\bar{M}_a$ and $\bar{\omega}$ are parameters of the inertia dynamics associated with an inertia magnitude and frequency threshold, in analogy to \eqref{set_point}.
For enhanced flexibility, a range of values was considered for $\bar{M}_a$ and $\bar{\omega}$, with aim to minimize the value of $\kappa$ where unstable behaviour is observed.
This setting lead to unstable behaviour for values of $\kappa$ as low as $2.8$, which suggests that in the considered simulation setting,  the conditions imposed by Assumption \ref{assum_inertia_growth} could have been relaxed by up to 2.8 times.
However, it should be noted that  Assumption \ref{assum_inertia_growth} imposes decentralized conditions, applicable to any network configuration, and hence some conservativeness is reasonable to be expected.

}

\subsubsection{Instability inducing single-bus varying inertia}
\begin{figure}[t!]
\includegraphics[scale=0.63]{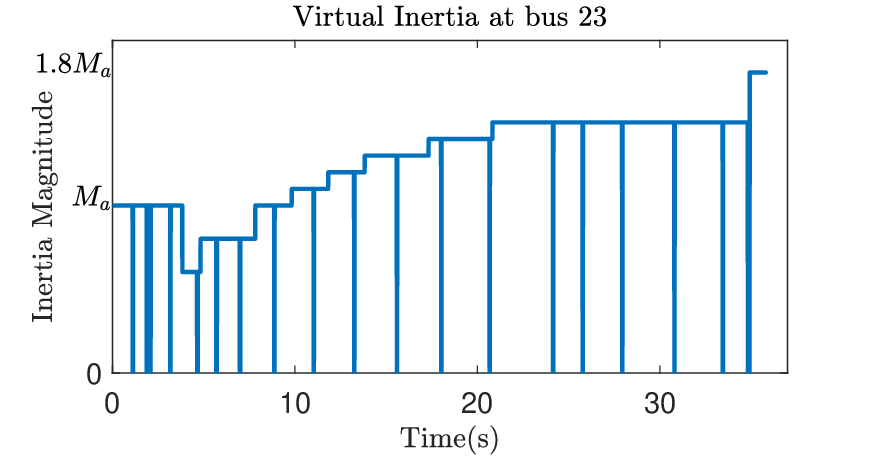}
\vspace{-5mm}
\caption{Virtual inertia at bus $23$ that resulted in unstable behaviour under constant inertia in the remaining power network. }
\label{Fig_single_inertia}
\end{figure}

\begin{figure}[t!]
\includegraphics[scale=0.62]{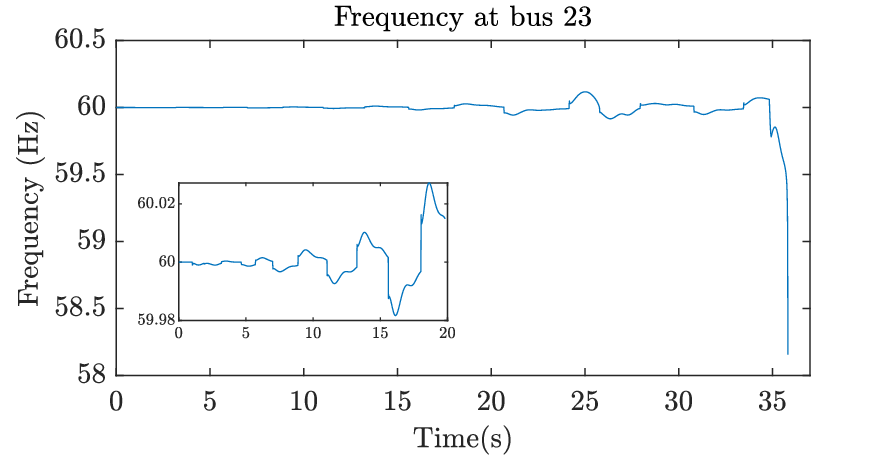}
\vspace{-5mm}
\caption{Frequency response at bus $23$ when varying virtual inertia is present at bus $23$ and no virtual inertia is present in the remaining buses.}
\label{Fig_single_instability}
\end{figure}

To demonstrate that local inertia variations may result in unstable behaviour, we performed simulations on the above described setting and considered a very small disturbance of magnitude $0.05$ p.u. (base 100MVA) at load buses $2$ and $3$  at $t=1$ second.
In addition, we considered no virtual inertia at all buses except bus $23$.
The design of the virtual inertia trajectory at bus $23$ followed the arguments in the proof of Theorem \ref{thm_instability}. 
In particular, the virtual inertia shifted between a set of large values and zero when the frequency deviations from the nominal value were large and small respectively.
The virtual inertia trajectory at  bus $23$ is depicted in Fig. \ref{Fig_single_inertia}.
As follows from  Fig. \ref{Fig_single_inertia}, the virtual inertia is piecewise-constant and takes zero values for short durations of time.
The frequency response for the considered case is presented in Fig. \ref{Fig_single_instability}, which depicts increasing frequency oscillations as time grows which eventually lead to instability after approximately 35 seconds.
The presence of increasing frequency deviations from equilibrium is in agreement with Assumption \ref{assum_instability}(ii), which is a main condition to deduce instability.
These results verify the analysis presented in Section \ref{Sec_Instability}.

\subsection{Simulations on the IEEE New York / New England 68-bus system}\label{subsection_NYNE}

The IEEE New York / New England system contains $52$ load buses serving different types
of loads including constant active and reactive loads and $16$ generation buses. The overall
system has a total real power of $16.41$ GW. 
For our simulation, we considered a step increase in demand of magnitude $2$ p.u. at load buses $2$ and $3$ at $t = 1$ second. The time precision of the simulation was $0.01$ seconds.

In analogy to Section \ref{subsection_NPCC}, we aimed to demonstrate how varying inertia may induce large frequency oscillations and that the application of the proposed conditions enables a stable response. We considered the following three cases:
\begin{enumerate}[(i)]
\item The presence of no virtual inertia.
\item The presence of virtual inertia in all generation buses, with dynamics described by \eqref{unstable_VI}, where $M_a = 20$ $s$.
\item The presence of virtual inertia in all generation buses with dynamics described by \eqref{Stable_VI}, where $M_a = 20$ $s$ and $\tau = 100$ $\text{s}^{-1}$.
\end{enumerate}
Case (ii) includes fast varying, frequency dependent inertia trajectories, which do not abide to the conditions presented in this paper. On the other hand, case (iii) includes inertia with variations bounded by the local strictness constant $\rho_j$, \an{in agreement with Assumption \ref{assum_inertia_growth}.}

\begin{figure}
\includegraphics[scale=0.63]{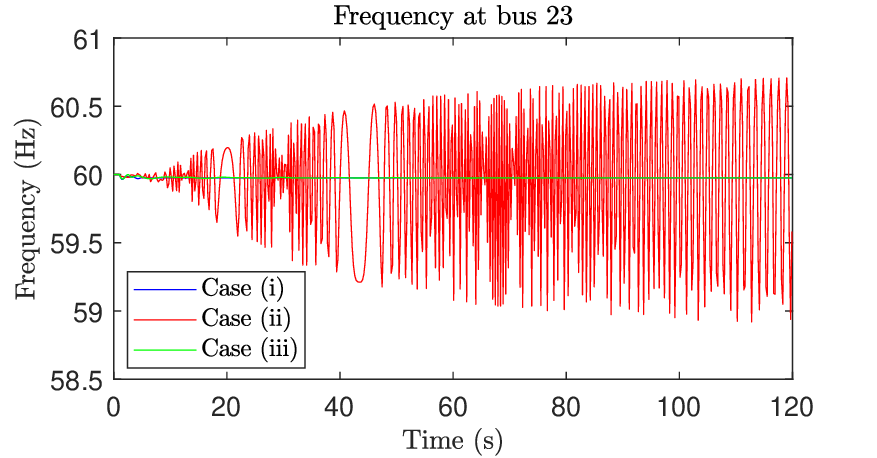}
\vspace{-5mm}
\caption{Frequency response at bus $23$ when: (i) no virtual inertia is present, (ii)~virtual inertia with dynamics described by \eqref{unstable_VI} is included in all generation buses and (iii) virtual inertia with dynamics described by \eqref{Stable_VI} is included in all generation buses.}
\label{Fig_frequency_16}
\vspace{-2mm}
\end{figure}

The frequency response at bus $23$ for the three considered cases is depicted in Fig. \ref{Fig_frequency_16}.
From Fig. \ref{Fig_frequency_16}, it follows that case (ii) yields an oscillatory frequency response, demonstrating how fast inertia oscillations may cause stability issues. On the other hand, a stable response was observed when case (iii) was implemented. 
The latter demonstrates how imposing a suitable bound on the rate of change of inertia may enable a stable response and verifies the analysis presented in Section \ref{Sec: Convergence}.

\begin{figure}[t!]
\includegraphics[scale=0.62]{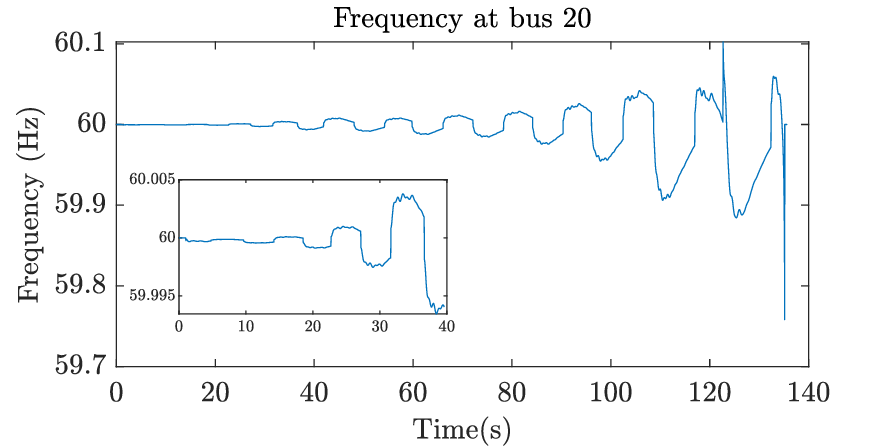}
\vspace{-5mm}
\caption{Frequency response at bus $20$ when varying virtual inertia is present at bus $57$ and no virtual inertia is present at the remaining buses.}
\label{Fig_single_instability_NYNE}
\vspace{-2mm}
\end{figure}

\begin{figure}
\centering
\includegraphics[scale=0.62]{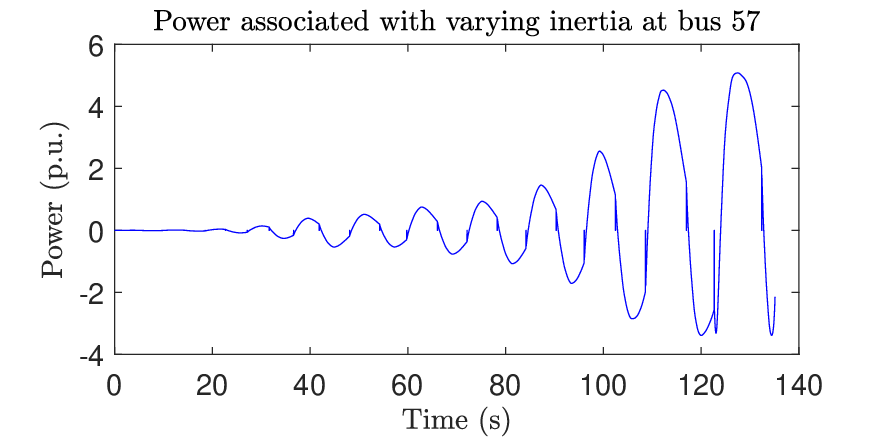}
\vspace{-5mm}
\caption{\ak{Power associated with varying inertia at bus $57$ 
that resulted in unstable behaviour under constant inertia in the remaining power network.}}
\label{Fig_NYNE_instability_power}
\end{figure}

To demonstrate how local inertia variations may result in unstable behaviour, we simulated the above described setting with a very small load disturbance of magnitude $0.01$ p.u. at load buses $2$ and $3$ at $t = 1$ second.
In addition, we considered varying inertia at bus $57$ and constant inertia at all remaining buses.
Similar to Section \ref{subsection_NPCC}, the design of the virtual inertia trajectory at bus $57$ followed the arguments in the proof of Theorem \ref{thm_instability}, i.e. large inertia values were employed during large frequency deviations and no virtual inertia was considered otherwise. 
The frequency response at a randomly selected bus (bus $20$) is presented in Fig. \ref{Fig_single_instability_NYNE}. 
Figure \ref{Fig_single_instability_NYNE} depicts growing oscillations, caused due to the inertia variations at bus $57$. 
Considering that such response follows from a normally negligible disturbance (of $2$ MW) demonstrates the capability of local inertia variations to cause instability in power networks and verifies the analysis presented in Section~\ref{Sec_Instability}.
\ak{Finally, the power associated with varying inertia at bus $57$ is depicted in Fig. \ref{Fig_NYNE_instability_power}. From Fig. \ref{Fig_NYNE_instability_power}, it follows that the power associated with varying inertia gradually increases in an oscillating manner, reaching an amount that is substantial, but still a small portion of the total network power, to yield unstable behaviour.
}


\section{Conclusion}\label{Sec: Conclusions}

We have investigated the stability properties of power networks with time-varying inertia within the primary frequency control timeframe.
In particular, we considered the interaction between the frequency dynamics and a wide class of nonlinear power supply dynamics at the presence of time-varying  inertia.
For the considered system,  we provided asymptotic stability guarantees under two proposed conditions.
The first condition required that the aggregate power supply dynamics at each bus satisfied a passivity related property.
The second condition set a constraint on the maximum rate of growth of inertia that was associated with the local power supply dynamics.
The proposed conditions are decentralized and applicable to arbitrary network topologies
and may be used for the design of practical guidelines for virtual inertia  schemes that will improve the reliability and enhance the stability properties of the power grid.
In addition, to demonstrate their applicability,  we explain how these conditions can be efficiently verified for linear power supply dynamics by solving a suitable linear matrix inequality optimization problem.
Our stability analysis is complemented with further analytic results that demonstrate how single-bus inertia variations may lead to instability.
\an{These results demonstrate that instability may arise from the interaction between frequency dynamics and inertia schemes, providing intuition on the underlying instability mechanisms, and highlighting the importance of incorporating the developed stability conditions on inertia scheme design.
} 
%
Numerical simulations on the  NPCC 140-bus and New York/ New England 68-bus systems offered additional motivation and validated our analytic results.
In particular, the simulation results   demonstrated how varying virtual inertia may induce large frequency oscillations and verified that the application of the proposed conditions resulted in a stable response.
In addition, they illustrate how single-bus inertia variations   may lead to unstable power system behaviour.

\section*{Appendix}

This appendix includes the proofs of Lemmas \ref{lemma_existence_uniqueness} and \ref{lemma_conv},  Theorems \ref{thm_conv} and \ref{thm_instability} and Corollary \ref{cor_instability}. Additionally, it includes Propositions \ref{prop_conv}, \ref{prop_finite_time} and \ref{prop_trajectories} that facilitate the proof of Theorem \ref{thm_instability}.

\emph{Proof of Lemma  \ref{lemma_existence_uniqueness}:}
Existence of a unique solution $\beta(t)$ to \eqref{sys1_compact}--\eqref{dynsys} for all $t \geq 0$ requires  \cite[Ch. 4.3]{haddad2011nonlinear} that
(i)  \eqref{sys1_compact}--\eqref{dynsys} is continuous in $t$ and $\beta$, and
(ii) \eqref{sys1_compact}--\eqref{dynsys} is locally Lipschitz in $\beta$ uniformly in $t \in [0, \infty)$.
Condition (i) is satisfied from Assumption \ref{assum_Virtual_Inertia} and the continuity of \eqref{sys1_compact}--\eqref{dynsys} in $\beta$.
Condition (ii) is satisfied since  \eqref{sys1_compact}--\eqref{dynsys} is locally Lipschitz in $\beta$ for given $t$ and the fact that $M_j(t) \geq M^0_{j} > 0$ for all times which allows a uniform local Lipschitz constant to be obtained for all $t \geq 0$.
Hence, for any initial condition $\beta(0) \in \mathbb{R}^m$ there exists a unique solution $\beta(t), t \geq 0$ to \eqref{sys1_compact}--\eqref{dynsys}.
 \hfill $\blacksquare$

\emph{Proof of Theorem \ref{thm_conv}:}
We will use Lyapunov arguments to prove Theorem \ref{thm_conv}.
First, we consider the  function
\begin{equation*}
V_F(M, \omega) =\frac{1}{2} \sum_{j \in \mathcal{N}} M_j(t) (\omega_j - \omega^*_j)^2,
\end{equation*}
with time derivative along trajectories of \eqref{sys1b} given by 
\begin{align}
&\dot{V}_F = \sum_{j \in \mathcal{N}} [\frac{\dot{M}_j}{2}  (\omega_j - \omega^*_j)^2 \nonumber \\ &+ (\omega_j - \omega^*_j) ( - p_j ^L + s_j
 - \sum_{k \in \mathcal{N}^s_j} p_{jk} + \sum_{l \in \mathcal{N}^p_j} p_{lj})],\label{VM_derivative}
\end{align}
In addition, we consider the function 
\begin{equation*}
V_{P}(\eta) = \sum_{(k,l) \in \mathcal{E}} \int_{\eta^*_{kl}}^{\eta_{kl}} B_{kl}(\sin \phi - \sin \eta^*_{kl}) d \phi.
\end{equation*} 
 From \eqref{sys1a} and \eqref{sys1c}, the derivative of $V_P$  is given by
\begin{align}\label{VP_derivative}
\dot{V}_P &= 
\sum_{(k,l) \in \mathcal{E}} B_{kl}(\sin \eta_{kl} - \sin \eta^*_{kl})(\omega_k - \omega_l) \nonumber\\ &= 
\sum_{(k,l) \in \mathcal{E}} (p_{kl} - p^*_{kl})(\omega_k - \omega_l).
\end{align}
Furthermore, from Assumption~\ref{assum2} and Definition \ref{Passivity_Definition}, it follows that there exist open neighbourhoods $\Omega_j$ of $\omega^*_j$ and $X_j$ of $x^{s,*}_j$ and  continuously differentiable, positive semidefinite functions $V_j(x^{s}_j)$  such that
\begin{multline}\label{Vd_derivative}
\hspace{0em} \dot{V}_j \leq ((-\omega_j) - (-\omega^*_j))(s_j - s^{*}_j) 
 - \rho_j((-\omega_j) - (-\omega^*_j))^2,
\end{multline}
where $\rho_j > 0$, for all $\omega_j \in \Omega_j, x^{s}_j \in X_j$ and for each $j \in \mathcal{N}$.
We now consider the following, \ank{positive-definite, continuously differentiable,} Lyapunov candidate function
\begin{equation}\label{Lyapunov_function}
V(M, \beta) = V_F(M, \omega) + V_{P}(\eta) + \sum_{j \in \mathcal{N}} V_j(x^{s}_j),
\end{equation}
reminding that $\beta = (\omega, \eta, x^{s})$.
Using \eqref{VM_derivative}--\eqref{Vd_derivative} it follows that when $\omega_j \in \Omega_j, x^{s}_j \in X_j, j \in \mathcal{N}$, then
\begin{equation}\label{dot_V}
\dot{V} \leq \sum_{j \in \mathcal{N}} (\dot{M}_j/2 - \rho_j) (\omega_j - \omega^*_j)^2 \leq  -W(\beta) \leq 0,
\end{equation}
where  the first inequality follows by applying \eqref{eqbr2} on \eqref{VM_derivative} and the second from Assumption \ref{assum_inertia_growth}, \ank{where $W(\beta) = \frac{1}{2}\sum_{j \in \mathcal{N}} \epsilon_j (\omega_j - \omega^*_j)^2$ is a continuous, positive semidefinite function.}
\ka{Moreover, we define $S := \{M \in R^{\mathcal{N}}: M^0_j \leq M_j \leq \bar{M}_j, j\in \mathcal{N}\}$ such that $M(t) \in S, t \geq 0$, where $\bar{M}_j$ is the upper bound of $M_j, j \in \mathcal{N}$ as follows from Assumption \ref{assum_Virtual_Inertia}.}

Function $V_F$ has a global minimum at $\omega = \omega^*$. In addition, Assumption \ref{assum1} guarantees the existence of a neighbourhood of $\eta^*$ where $V_P$ is increasing, which suggests that $V_P$ has a strict local minimum at $\eta^*$. 
Moreover, from Assumption \ref{assum2} and Definition \ref{Passivity_Definition} it follows that each $V_j, j \in \mathcal{N}$ has a strict local minimum at $x^{s,*}_j$. Hence, $V$ has a  local minimum at $\beta^* = (\omega^*, \eta^*, x^{s,*})$ that is independent of the value of $M(t)$, since $M_j(t) \geq M^0_j > 0, t \geq 0, j \in \mathcal{N}$. In addition, $\beta^*$ is a strict minimum associated with the states, i.e. locally the set $\{\bar{\beta} : \bar{\beta} = \argmin_{\beta} V(M, \beta), \forall M \in S\}$ contains $\beta^*$ only. We can now choose a neighbourhood of $\beta^*$, denoted by $\bar{B}$, such that 
(i) $\omega_j \in \Omega_j, j \in \mathcal{N}$,
(ii) $x^s_j \in X_j, j \in \mathcal{N}$, and
(iii) all $x^s_j, j \in \mathcal{N}$ lie in their respective neighbourhoods $\Psi_j$, as defined in Section \ref{Sec: Generation and demand dynamics}.

\ank{Then, we invoke \cite[Theorem 4.1 (i)]{haddad2011nonlinear}  to deduce Lyapunov stability for \eqref{sys1_compact}--\eqref{dynsys} around $\beta^*$.
To show this, we define, without loss of generality, \(\tilde{V}(M, \beta - \beta^*) = V(M, \beta)\), and note that it satisfies (i) \(\tilde{V}(M, 0) = 0\) for \(M \in S\), (ii) inequality \eqref{dot_V}, and (iii)~\(\alpha(\|\beta - \beta^*\|) \leq \tilde{V}(M, \beta - \beta^*)\) for some class \(\mathcal{K}\) function\footnote{A  continuous function $\alpha: [0, \tilde{\alpha}] \rightarrow [0, \infty)$, where $\tilde{\alpha}$ is a positive constant, is said to belong to class $\mathcal{K}$ if it is strictly increasing and satisfies $\alpha(0) = 0$.} \(\alpha\). Such class \(\mathcal{K}\) function can be constructed by considering \(\tilde{V}(M, \beta - \beta^*)\) and fixing \(M = M^0\), and leveraging that $\tilde{V}$ is continuously differentiable with a strict local minimum at \(\beta = \beta^*\).
This suggests  that for every  \(\varepsilon > 0\) such that $\Lambda = \{ \beta: \|\beta - \beta^*\| \leq \epsilon \} \subseteq \bar{B}$, 
there exists \(\delta = \delta(\varepsilon) > 0\) such that \(\beta(0) \in \tilde{\Xi} = \{ \beta: \|\beta - \beta^*\| \leq \delta \} \) implies that \(\beta(t) \in \Lambda\) for all \( t \geq 0 \) (see \cite[Definition 4.1]{haddad2011nonlinear}).
\ka{Note that $\delta$ above can be defined as a function of $\varepsilon$ only, since $M$ is bounded, as follows from Assumption \ref{assum_Virtual_Inertia}.} 
For $(M, \beta) \in S \times \Lambda$, $V$ is a non-increasing function with a strict local minimum associated with the states at $\beta^*$.
 Hence, $\lim_{t \rightarrow \infty} V (M(t), \beta(t))$ exists and is finite \cite[Theorem 2.10]{haddad2011nonlinear}.
This suggests that $\lim_{t \rightarrow \infty} \int_0^t W (\beta(\tau)) d\tau$ also exists and is finite, since $\int_0^t W(\beta(\tau)) d\tau \leq - \int_0^t \dot{V} (M(\tau), \beta(\tau)) d\tau = V(M(0), \beta(0)) - V(M(t), \beta(t))$.
Then, we note that within the compact set $\Lambda$, for any $M(t) \in S, t \geq 0$, the dynamics \eqref{sys1_compact}--\eqref{dynsys} are Lipschitz continuous uniformly in $M$.
Then, it follows that within $\Lambda$, $\| \beta(t_2) - \beta(t_1) \| \leq \tilde{L} \varepsilon (t_2 - t_1), t_2 \geq t_1$, which allows to deduce that $\beta$ is uniformly continuous since for any 
$\tilde{\epsilon}$, 
we can define $\tilde{\delta} = \tilde{\epsilon}/(\tilde{L} \epsilon)$ such that 
$\| \beta(t_2) - \beta(t_1) \| <~\tilde{\epsilon},$ when $t_2 - t_1 < \tilde{\delta}$.
Hence, since $W$ is continuous, we can deduce that within $\Lambda$, $W(\beta(t))$ is uniformly continuous for every $t \geq 0$.
}

\ank{Then, using Barbalat's Lemma \cite[Lemma 4.1]{haddad2011nonlinear}, it follows that $W(\beta(t)) \rightarrow 0$ as $t \rightarrow \infty$.
Therefore, we deduce that  for any $(M(0), \beta(0)) \in S \times \tilde{\Xi}$ it follows that $\beta \rightarrow Q$ as $t \rightarrow \infty$, 
where $Q = \Lambda \cap \{\beta: W(\beta) = 0 \}$.}
Within this set, it holds that $\omega_j = \omega^*_j, j \in \mathcal{N}$ from \eqref{dot_V}. In addition, the definitions in Section \ref{Sec: Generation and demand dynamics} suggest that $\omega = \omega^*$ implies the convergence of $x^s$ to $x^{s,*}$.

\ka{Next, we will use Barbalat's Lemma \cite[Lemma 4.1]{haddad2011nonlinear} again to show that 
$\dot{\omega}_j \rightarrow 0$ as $t \rightarrow \infty$ for $j \in \mathcal{N}$.
First, we show that  $\dot{\omega}_j(t)$ is uniformly continuous for all $j \in \mathcal{N}$,
by differentiating \eqref{sys1b} and using \eqref{ssys} to obtain
\begin{multline}\label{ddot_omega}
\ddot{\omega}_j =  
(\dot{s}_j - \dot{M}_j \dot{\omega}_j  - \sum_{k \in \mathcal{N}^s_j} B_{jk}\cos(\eta_{jk}) \dot{\eta}_{jk} \\
+  \sum_{l \in \mathcal{N}^p_j} B_{lj}\cos(\eta_{lj}) \dot{\eta}_{lj}) / M_j(t) \\
\leq
|\dot{s}_j - \dot{M}_j \dot{\omega}_j  - \sum_{k \in \mathcal{N}^s_j} B_{jk}\cos(\eta_{jk}) \dot{\eta}_{jk} \\
+  \sum_{l \in \mathcal{N}^p_j} B_{lj}cos(\eta_{lj}) \dot{\eta}_{lj}| / M^0_j, 
j\in \mathcal{N},
\end{multline}
where the inequality follows by noting that $M_j(t) \geq M^0_j \geq 0$.
Moreover, note that the magnitude of $\dot{M}_j$ is bounded, as follows from Assumption \ref{assum_Virtual_Inertia} and for $\beta \in \Lambda$, 
(a) the magnitudes of $\dot{\eta}_{jk}, (j,k) \in \mathcal{E}^s_j$, and $\dot{\omega}_j, j \in \mathcal{N}$ are bounded, as follows from \eqref{sys1a_compact}, \eqref{sys1b_compact}, and $M_j(t) \geq M^0_j \geq 0$, 
and (b) the magnitude of $\dot{s}_j$ is bounded, as follows by noting that $\dot{s}_j = \frac{dg_j(x^s_j,-\omega_j)}{dx^s_j}f_j(x^s_j,-\omega_j) + \frac{dg_j(x^s_j,-\omega_j)}{d\omega_j}\dot{\omega}_j$ and that  $g_j$ is locally Lipschitz\footnote{Although $g_j$ is only assumed to be locally Lipschitz (so its partial derivatives do not necessarily exist everywhere), local Lipschitzness implies differentiability almost everywhere (Rademacher’s theorem). 
Since $(x^s_j(t),-\omega_j(t))$ are Lipschitz in time when the trajectories lie within the compact set $\Lambda$, $s_j(t)=g_j(x^s_j(t),-\omega_j(t))$ admits a derivative almost everywhere that satisfies the chain rule at points of differentiability. 
In particular, within $\Lambda$,  $|\dot s_j(t)|$ is essentially bounded, which is sufficient for the uniform-continuity arguments in the proof.}.
Hence, for solutions that satisfy $\beta \in \Lambda$, it follows that 
$\|\dot{\omega}_j(t_2) - \dot{\omega}_j(t_1)\| = \|
\int_{t_1}^{t_2} \ddot{\omega}_j(\tau) d\tau\| \leq \bar{L}_j (t_2 - t_1), t_2 \geq t_1$, where $\bar{L}_j$ is the bound on the magnitude of $\ddot{\omega}_j$.
Therefore, for any $\bar{\epsilon}$ there exists $\bar{\delta} = \bar{\epsilon}/\bar{L}_j$ such that  $\|\dot{\omega}_j(t_2) - \dot{\omega}_j(t_1)\| < \bar{\epsilon}$ when $\|t_2 - t_1\| < \bar{\delta}$.
Hence, $\dot{\omega}_j(t)$ is uniformly continuous. Moreover, the convergence of $\omega_j$ to a constant value suggests that $\lim_{t \rightarrow \infty} \int_{0}^{t}\dot{\omega}_j (\tau) d\tau$ exists and is finite.
Hence, we can invoke Barbalat's Lemma to deduce that
$\lim_{t \rightarrow \infty} \dot{\omega}_j (t) = 0, j \in \mathcal{N}$. Moreover, noting that $M(t)$ is bounded, as follows from Assumption \ref{assum_Virtual_Inertia}, it follows that $\lim_{t \rightarrow \infty} M_j(t) \dot{\omega}_j (t) = 0, j \in \mathcal{N}$.
}


\ka{The above arguments suggest the convergence of $\beta$ along trajectories of \eqref{sys1_compact}--\eqref{dynsys} to a subset of $Q$, that we denote by $\Theta$, where $\omega_j = \omega^*_j, j \in \mathcal{N}$, noting also that $\omega^*_i = \omega^*_j, (i,j) \in \mathcal{E}$, as follows from Definition \ref{eqbrdef}, $x^s = x^{s,*}$, and \eqref{eqbr2} hold, where the latter follows from $\lim_{t \rightarrow \infty} M_j(t) \dot{\omega}_j (t) = 0, j \in \mathcal{N}$.} 
\ka{The above suggest that all points  within $\Theta$ also satisfy  \eqref{eqbr1} and \eqref{eqbr3}.}
\an{Hence, all equilibrium conditions \eqref{eqbr} hold within $\Theta$.}
Therefore, we conclude  that 
there exists a set $\tilde{\Xi}$ containing $\beta^*$ such that solutions $\beta(t), t \geq 0$ to \eqref{sys1_compact}--\eqref{dynsys} initiated in $\tilde{\Xi}$ converge to the set of equilibrium points of \eqref{sys1_compact}--\eqref{dynsys}. 
Letting $\Xi$ be any open neighbourhood of $\beta^*$ within $\tilde{\Xi}$ completes the proof.
 \hfill $\blacksquare$

 \emph{Proof of Lemma \ref{lemma_conv}:} 
 The proof follows directly from the proof of Theorem \ref{thm_conv} by letting $V_{P}(\eta) =  \sum_{(k,l) \in \mathcal{E}} B_{kl}(\eta_{kl} - \eta^*_{kl})^2/2$. In particular, it follows that locally around the considered equilibrium, solutions to \eqref{dynsys}, \eqref{sysl_compact} satisfy \eqref{dot_V}. The rest of the arguments follow in analogy to the proof of Theorem \ref{thm_conv}.
\hfill $\blacksquare$

%


To facilitate the analysis associated with the instability results presented in Section \ref{Sec_Instability}, we consider the effect of having some bus $k$ with a fixed frequency $\bar{\omega}$ within system \eqref{dynsys}, \eqref{sysl_compact}. 
The dynamics of such a system are described below:
\begin{subequations} \label{sys2_compact}
\begin{align}
\dot{\eta} &= H^T \omega, \label{sys2a_compact}
\\
 M_j \dot{\omega}_j&= \hspace{-0.5mm} - p^L_j \hspace{-0.5mm}+\hspace{-0.5mm} s_j \hspace{-0.5mm}-\hspace{-1mm} \sum_{k \in \mathcal{N}^s_j}\hspace{-0.5mm} p_{jk} +\hspace{-1mm} \sum_{l \in \mathcal{N}^p_j}\hspace{-0.5mm} p_{lj}, j\in \mathcal{N}\setminus \{k\}, \label{sys2b_compact}
\\
\omega_k &= \bar{\omega}, \label{sys2d_compact} \\
p &= B    \eta. \label{sys2c_compact}
\end{align}
\end{subequations}
The system represented by \eqref{sys2_compact} aims to describe the behaviour of \eqref{sysl_compact}, when the inertia of bus $k$ is infinite and $\omega_k (0) = \bar{\omega}$.
Although the assumption of infinite inertia is unrealistic, the trajectories of \eqref{dynsys}, \eqref{sys2_compact} approximate those of \eqref{dynsys},  \eqref{sysl_compact}  for a long time interval when the inertia at bus $k$ is sufficiently large. The latter enables to explore several properties of  \eqref{dynsys},  \eqref{sysl_compact} and facilitates our instability analysis.

It should be noted that there is a direct relation between \eqref{sys2_compact} and the analysis presented in Section \ref{Sec_Instability}, since it can be trivially shown that $\gamma$-points, defined in Definition \ref{eqbrdef_2}, coincide with the equilibria to \eqref{dynsys}, \eqref{sys2_compact}.

The following proposition  demonstrates the convergence of solutions to \eqref{dynsys}, \eqref{sys2_compact} to the set of its equilibria, similarly to Theorem \ref{thm_conv}. In addition, it provides conditions that allow to deduce convergence to an equilibrium point of \eqref{dynsys}, \eqref{sys2_compact}.
It should be clarified that within Proposition \ref{prop_conv}, and also Proposition \ref{prop_finite_time} below, Assumptions \ref{assum_Virtual_Inertia}, \ref{assum_inertia_growth} refer to all buses besides bus $k$. The latter follows since no inertia is defined for  bus $k$.

\begin{proposition}\label{prop_conv}
Let Assumptions \ref{assum_Virtual_Inertia} and \ref{assum_inertia_growth} hold and consider an equilibrium of  \eqref{dynsys}, \eqref{sys2_compact}, \ank{denoted by $\beta^*$}, where Assumption \ref{assum2} holds.
Then, 
\begin{enumerate}[(i)]
\item 
 {there exists an open neighbourhood $\Xi$ containing that equilibrium} such that solutions $\beta(t), t \geq 0$ \eqref{dynsys}, \eqref{sys2_compact} initiated in $\Xi$ asymptotically converge to the set of equilibria of \ank{\eqref{dynsys}, \eqref{sys2_compact}},
\item   \ank{for every $\epsilon > 0$}
 {there exists an open neighbourhood $\Xi$ containing that equilibrium such that solutions $\beta(t), t \geq 0,$ to \eqref{dynsys}, \eqref{sys2_compact} initiated in $\Xi$ satisfy $\beta(t) \in \Lambda = \{\beta: \|\beta - \beta^*\| \leq \epsilon\}$.}
Moreover,  if Assumption \ref{assum2} holds for all equilibria within $\Lambda$, then 
solutions $\beta(t), t \geq 0,$ to \eqref{dynsys}, \eqref{sys2_compact} initiated in $\Xi$ converge to an equilibrium point within $\Lambda$.
\end{enumerate} 
\end{proposition}

\emph{Proof of Proposition \ref{prop_conv}:}
The proof is split in two parts, regarding each statement in Proposition \ref{prop_conv}. 

\emph{Part (i):}  The proof follows by using Lyapunov  arguments, similar to the proof of Theorem \ref{thm_conv}.
First, we consider the  function
$
V_G(M, \omega) = \frac{1}{2} \sum_{j \in \mathcal{N} \setminus \{k\}} M_j(t) (\omega_j - \omega^*_j)^2,
$
with time derivative along trajectories of \eqref{sys2b_compact} given by 
\begin{align}
&\dot{V}_G = \sum_{j \in \mathcal{N} \setminus \{k\}} [\frac{\dot{M}_j}{2}  (\omega_j - \omega^*_j)^2 \nonumber \\ &+ (\omega_j - \omega^*_j) ( - p_j ^L + s_j
 - \sum_{k \in \mathcal{N}^s_j} p_{jk} + \sum_{l \in \mathcal{N}^p_j} p_{lj})].\label{VG_derivative}
\end{align}
In addition, we consider the function $V_{H}(\eta) = \sum_{(k,l) \in \mathcal{E}} B_{kl}(\eta_{kl} - \eta^*_{kl})^2/2$, with derivative given by \eqref{sys2a_compact} and \eqref{sys2c_compact} as 
\begin{align}\label{VH_derivative}
\dot{V}_H &= 
\sum_{(m,l) \in \mathcal{E}} (p_{ml} - p^*_{kl})(\omega_m - \omega_l).
\end{align}
Furthermore, from Assumption~\ref{assum2} and Definition \ref{Passivity_Definition}, it follows that there exist open neighbourhoods $\Omega_j$ of $\omega^*_j$ and $X_j$ of $x^{s,*}_j$ and  continuously differentiable, positive semidefinite functions $V_j(x^{s}_j)$  such that \eqref{Vd_derivative} holds.

We now consider the following Lyapunov candidate 
\begin{equation*}
V(M, \beta) = V_G(M, \omega) + V_{H}(\eta) + \sum_{j \in \mathcal{N}} V_j(x^{s}_j),
\end{equation*}
reminding that $\beta = (\omega, \eta, x^{s})$.
Using \eqref{Vd_derivative}, \eqref{VG_derivative}, \eqref{VH_derivative} it follows that 
\begin{equation}\label{dot_V2}
\dot{V} \leq \sum_{j \in \mathcal{N} \setminus \{k\}} (\dot{M}_j/2 - \rho_j) (\omega_j - \omega^*_j)^2 \leq 0,
\end{equation}
by applying \eqref{eqbr2} on \eqref{VG_derivative} and using Assumption \ref{assum_inertia_growth}. 

Then, similar arguments as in the proof of Theorem \ref{thm_conv}, it follows that solutions initiated in \ank{$\Xi$ converge to the set of equilibria of \eqref{dynsys}, \eqref{sys2_compact}}.


\emph{Part (ii):} 
\ank{Using similar arguments as in the proof of Theorem \ref{thm_conv}, it follows that the considered equilibrium point is Lyapunov stable \cite[Definition 4.1]{haddad2011nonlinear}. Hence, for every $\epsilon > 0$ such that $\Lambda = \{\beta: \|\beta - \beta^*\| \leq \epsilon\} \subseteq \bar{B}$
  there exists a set $\tilde{\Xi}$ containing that equilibrium such that solutions $\beta(t), t \geq 0,$ to \eqref{dynsys}, \eqref{sys2_compact} initiated in $\tilde{\Xi}$ satisfy $\beta(t) \in \Lambda$.}
Now if Assumption \ref{assum2} holds for all equilibria within $\Lambda$ it follows that these equilibria are also Lyapunov stable.
The latter allows to use \cite[Th. 4.20]{haddad2011nonlinear} to deduce \ank{the existence of an open subset $\Xi$ within $\tilde{\Xi}$ containing $\beta^*$ such that solutions initiated in $\Xi$} converge to an equilibrium point within $\Lambda$.
%
 \hfill $\blacksquare$

The following proposition, shows that  system \eqref{dynsys}, \eqref{sys2_compact} has arbitrarily long periods of time where all bus frequencies lie within a ball of size $\epsilon$ from $\bar{\omega}$, for any positive value of $\epsilon$.

\begin{proposition}\label{prop_finite_time}
Let Assumptions \ref{assum_Virtual_Inertia} and \ref{assum_inertia_growth} hold and consider an equilibrium of  \eqref{dynsys}, \eqref{sys2_compact} where Assumption \ref{assum2} holds.
Then, 
there exists an open neighbourhood $\Xi$ containing that equilibrium such that for any solution $\beta(t), t \geq 0$ to \eqref{dynsys}, \eqref{sys2_compact} initiated in $\Xi$
and any $\tau, \epsilon \in \mathbb{R}_+$ there exists some $\hat{\tau} \in \mathbb{R}_+$
such that $\omega_j(t) \in \mathcal{B}(\bar{\omega}, \epsilon), t \in [\hat{\tau}, \hat{\tau} + \tau], j \in \mathcal{N}$.
\end{proposition} 
 
 \emph{Proof of Proposition \ref{prop_finite_time}:}
 First note that the results presented in Proposition \ref{prop_conv} hold, since all associated assumptions are  satisfied.
Then, from \eqref{dot_V2} and the Lipschitz continuity of \eqref{dynsys}, \eqref{sys2_compact}, it follows that  for any $\tau, \epsilon \in \mathbb{R}_+$  there exists some $\delta > 0$ such that for any $\bar{\tau}$, if $\omega_j(t) \notin \mathcal{B}(\bar{\omega}, \epsilon)$ for some $j \in \mathcal{N}$ and some $t \in  [\bar{\tau}, \bar{\tau} + \tau]$ then $V(\bar{\tau} + \tau) \leq V(\bar{\tau}) - \delta$.
Noting that $V(0)$ is bounded and $V(t) \geq 0, t \geq 0$ allows to deduce the above result by contradiction.
  \hfill $\blacksquare$

The following proposition states the existence of a virtual inertia trajectory $M^v_k(t), t \geq 0$ such that solutions to \eqref{dynsys}, \eqref{sysl_compact} and \eqref{dynsys}, \eqref{sys2_compact} are arbitrarily close for arbitrarily long, but finite, time intervals.
For convenience in presentation, we consider $\omega_k$ as a state of system  \eqref{dynsys}, \eqref{sys2_compact} that keeps a constant value at all times. The latter suggests that the statement regarding the initial conditions of \eqref{dynsys}, \eqref{sysl_compact} and \eqref{dynsys}, \eqref{sys2_compact} implies that $\omega_k(0)  = \bar{\omega}$. In addition, by slightly abusing notation, we denote trajectories of \eqref{dynsys}, \eqref{sysl_compact} and \eqref{dynsys}, \eqref{sys2_compact} by $\beta(t)$ and $\hat{\beta}(t)$ respectively.

\begin{proposition}\label{prop_trajectories}
Let Assumption \ref{assum_Virtual_Inertia} hold. Then, for any $\tau, \epsilon \in~\mathbb{R}_+$ and any trajectory for $M^v_j(t), j \in \mathcal{N} \setminus \{k\}, t \geq 0$, there exists a trajectory for $M^v_k(t), t \geq 0$ such that  solutions to  \eqref{dynsys}, \eqref{sysl_compact} and \eqref{dynsys}, \eqref{sys2_compact} with $\beta(0) = \hat{\beta}(0)$ satisfy $\norm{\beta(t) - \hat{\beta}(t)}\leq \epsilon, t \in [0, \tau]$.
\end{proposition}

\emph{Proof of Proposition \ref{prop_trajectories}:} 
First note that both  \eqref{dynsys}, \eqref{sysl_compact} and \eqref{dynsys}, \eqref{sys2_compact} are locally Lipschitz due to Assumption \ref{assum_Virtual_Inertia} and the conditions on \eqref{dynsys}. 
The proof follows by considering \eqref{sysl_compact} and noting that when $M^v_k \rightarrow \infty$ the frequency derivative at bus $k$ tends to zero and hence the frequency at bus $k$ is constant.
Hence, for any bounded trajectory for $M^v_j(t), j \in \mathcal{N} \setminus \{k\}, t \geq 0$, when  $\omega_k(0) = \bar{\omega}$  and $M^v_k = \infty$ the dynamics described by   \eqref{sysl_compact} and \eqref{sys2_compact} are identical.
The latter suggests the existence of sufficiently large, but finite, value $\bar{M}$ such that for any finite time $\tau$ and any $\epsilon > 0$, $M^v_k(t) \geq  \bar{M}, t \in [0, \tau]$ implies that $\norm{\beta(t) - \hat{\beta}(t)}\leq \epsilon, t \leq \tau$. The latter completes the proof.
 \hfill $\blacksquare$

\emph{Proof of Theorem \ref{thm_instability}:}
To prove Theorem \ref{thm_instability}, we will define an iterative process and provide properties for the  trajectory of $M_k$  such that there exist sequences of time instants $\hat{t}_i, i \in \mathbb{N}_+$ and positive values $\phi_i, i \in \mathbb{N}_+$ satisfying $\hat{t}_i > \hat{t}_j, i > j$ and $\phi_i  > \phi_j, i > j$ such that $\beta(\hat{t}_i) \notin \mathcal{B}(\gamma, \phi_i), \gamma \in \an{\Gamma_k}(\omega^{s,\ast}), i \in \mathbb{N}_+$. 
In addition, we will demonstrate the existence of some finite iteration $n$ such that $\phi_n \geq \Phi$, where $\Phi > 0$ is defined in Assumption \ref{assum_instability}.



From the theorem statement, it is assumed that $|\omega_k(0) - \omega^{s,\ast}| = \delta$. Then consider Proposition \ref{prop_trajectories}, \an{letting  $\bar{\omega} = \omega_k(0)$}. The latter claims that for any finite $\epsilon_1, \tau$, there exists a trajectory for $M^v_k(t)$ such that solutions to \eqref{dynsys}, \eqref{sysl_compact} and \eqref{dynsys}, \eqref{sys2_compact} with the same initial conditions have a distance of at most $\epsilon_1$, for $t \leq \tau$.
In addition,   from Proposition \ref{prop_finite_time} it follows that for any $\tau_2, \epsilon_2 \in \mathbb{R}_+$ there exists some $\hat{\tau}_2 \in \mathbb{R}_+$
such that $\omega_j(t) \in \mathcal{B}(\bar{\omega}, \epsilon_2), t \in [\hat{\tau}_2, \hat{\tau}_2 + \tau_2], j \in \mathcal{N}$.
%
Assumption \ref{assum_instability}(i) suggests the existence of some finite $\hat{\tau}$ such that the previous statement implies \an{for \eqref{dynsys}, \eqref{sys2_compact} that
 $\beta(\hat{\tau}) \in  \mathcal{B}(\gamma, {\epsilon'}), \gamma \in \an{\Gamma_k}(\bar{\omega}),$} $\an{{\epsilon}' > \epsilon_2}$.
  Using the previous arguments 
  suggests that solutions to \eqref{dynsys}, \eqref{sysl_compact} \an{satisfy 
$\beta(\hat{\tau}) \in  \mathcal{B}(\gamma, {\hat{\epsilon}}), \gamma \in \Gamma_k(\bar{\omega}), \hat{\epsilon} = \epsilon' + \epsilon_1$.}  
The latter additionally requires that $\tau \geq \hat{\tau}$ which can be achieved by suitably selecting $M^v_k$ following Proposition~\ref{prop_trajectories}.

 We now define $\hat{t}_1 =  \tau$ and $\bar{\omega}^1 = \delta$. From Assumption \ref{assum_instability}(ii), there exists some $\tau^1$ such that $M_k(t) = M^0_k, t \in [\hat{t}_1, \hat{t}_1 + \tau^1]$ implies that solutions to \eqref{dynsys}, \eqref{sysl_compact} satisfy $|\omega_k (\hat{t}_1 + \tau^1) - \omega^{s,\ast}| > \bar{\omega}^1 + \bar{\epsilon}, \bar{\epsilon} > \hat{\epsilon}$.
 If for some $t \in [\hat{t}_1, \hat{t}_1 + \tau^1]$ it holds that $\beta(\tilde{t}) \notin \mathcal{B}(\gamma, \Phi),  \gamma \in \an{\Gamma_k}(\omega^{s,\ast})$ then the proof is complete.
Otherwise, we let $\phi_1 = \bar{\omega}^1 + \bar{\epsilon}$, set a sufficiently large value for $M^v_k$ at $t = \hat{t}_1 + \tau^1$, as follows from Proposition \ref{prop_trajectories},  and repeat the above process iteratively\footnote{It should be noted that Assumption \ref{assum_Virtual_Inertia} requires that inertia trajectories are globally Lipschitz in time. Assumption \ref{assum_Virtual_Inertia} is satisfied by considering trajectories where the virtual inertia linearly changes between the zero and considered sufficiently large values and vice versa within some time duration $\bar{\delta}$ prior to $\hat{t}_i$ and $\hat{t}_i + \tau^i$ respectively, for some suitably selected value for  $\bar{\delta}$, noting that no bound is assumed on the Lipschitz constant.
  }, 
  using Assumption \ref{assum_instability}(iii) to deduce the convergence arguments\footnote{In particular, Assumption \ref{assum_instability}(iii) suggests that the regions where the asymptotic stability and passivity arguments used  to define the set $\Lambda$ in the proof of Proposition \ref{prop_conv} hold, are supersets of the regions where the associated trajectories considered in the proof arguments lie. } associated with Proposition \ref{prop_conv}.
This process creates a sequence of time instants  $\hat{t}_l$,  associated with each iteration $l \in \mathbb{N}_+$, satisfying $\hat{t}_{l+1} > \hat{t}_l$,  such that $|\omega_k(t_{l+1}) - \omega^{s,\ast}| \geq |\omega_k(t_l) - \omega^{s,\ast}| + \bar{\epsilon} - \hat{\epsilon}$, which implies that $\phi_{l+1}  > \phi_l  + \bar{\epsilon} - \hat{\epsilon}$.

Since the value of $\Phi$ is bounded, this will be reached in a finite amount of iterations (no more than $\Phi / (\bar{\epsilon} - \hat{\epsilon})$). Since the time required for each iteration is finite, then there exists some finite time  $\bar{t}$ such that 
$\beta(\bar{t}) \notin \mathcal{B}(\gamma, \Phi),  \gamma \in \an{\Gamma_k}(\omega^{s,\ast})$.
Noting that no assumption is made for the magnitude of $\delta$, and hence that the above arguments hold for any $\delta > 0$,
  completes the proof.
 \hfill $\blacksquare$

\emph{Proof of Corollary \ref{cor_instability}:}
The proof follows directly from Theorem \ref{thm_instability} which makes the same assumptions.
 In particular, an unstable equilibrium is defined as an equilibrium that is not stable, e.g. \cite[Dfn. 4.1]{haddad2011nonlinear}. 
For a system with state $x$, 
a stable equilibrium satisfies the property that for any $\epsilon > 0$, there exists  some $\delta > 0$ such that  $\norm{x(0)} < \delta$ implies $\norm{x(t)} < \epsilon, t \geq 0$.
Theorem \ref{thm_instability} states that for any $\delta > 0$, there exists an inertia trajectory at bus $k$ such that $|\omega_k - \omega^{s,\ast}| \geq \Phi$. The latter suggests that under specific inertia trajectories, there  exists some $\epsilon$ (i.e. any $\epsilon < \Phi$) such that there does not exist any $\delta > 0$ such that the resulting trajectories for \eqref{dynsys}, \eqref{sysl_compact} are bounded by $\epsilon$.
The latter completes the proof. 
 \hfill $\blacksquare$

\ank{
\section*{Acknowledgement}
The authors express their gratitude to the anonymous reviewers for their valuable comments and suggestions, which helped to improve the quality of the paper.
}

\balance
\bibliography{andreas_bib}

\end{document}